\documentclass[10pt,twoside,a4paper]{article}

\usepackage[a4paper,centering,body={31pc,49pc}]{geometry}

\usepackage{multicol}
\usepackage{amsmath}
\usepackage{amsfonts}
\usepackage{mathtools}
\usepackage{bm}
\usepackage{fancyhdr}
\usepackage{hyperref}

\usepackage{stmaryrd}

\usepackage{mathptmx}

\renewcommand{\sectionmark}[1]%
 {\markboth{\thesection:\ #1}{}}

\renewcommand{\Re}{\mathrm{Re}}%
\providecommand{\Ent}[1]{\lfloor #1 \rfloor}

\numberwithin{equation}{section}

\begin{document}

\thispagestyle{empty}

\title{\bfseries On the Analyticity of Laguerre Series}

\author{Ernst Joachim Weniger \\
  Institut f\"ur Physikalische und Theoretische Chemie \\
  Universit\"at Regensburg, D-93040 Regensburg, Germany}

\date{To appear in Journal of Physics A \\
  Accepted: 27 August 2008}

\maketitle

\begin{abstract}
  The transformation of a Laguerre series $f (z) = \sum_{n=0}^{\infty}
  \lambda_{n}^{(\alpha)} L_{n}^{(\alpha)} (z)$ to a power series $f (z) =
  \sum_{n=0}^{\infty} \gamma_{n} z^{n}$ is discussed. Since many
  nonanalytic functions can be expanded in terms of generalized Laguerre
  polynomials, success is not guaranteed and such a transformation can
  easily lead to a mathematically meaningless expansion containing power
  series coefficients that are infinite in magnitude. Simple sufficient
  conditions based on the decay rates and sign patters of the Laguerre
  series coefficients $\lambda_{n}^{(\alpha)}$ as $n \to \infty$ can be
  formulated which guarantee that the resulting power series represents
  an analytic function. The transformation produces a mathematically
  meaningful result if the coefficients $\lambda_{n}^{(\alpha)}$ either
  decay exponentially or factorially as $n \to \infty$. The situation is
  much more complicated -- but also much more interesting -- if the
  $\lambda_{n}^{(\alpha)}$ decay only algebraically as $n \to \infty$. If
  the $\lambda_{n}^{(\alpha)}$ ultimately have the same sign, the series
  expansions for the power series coefficients diverge, and the
  corresponding function is not analytic at the origin. If the
  $\lambda_{n}^{(\alpha)}$ ultimately have strictly alternating signs,
  the series expansions for the power series coefficients still diverge,
  but are summable to something finite, and the resulting power series
  represents an analytic function. If algebraically decaying and
  ultimately alternating Laguerre series coefficients
  $\lambda_{n}^{(\alpha)}$ possess sufficiently simple explicit
  analytical expressions, the summation of the divergent series for the
  power series coefficients can often be accomplished with the help of
  analytic continuation formulas for hypergeometric series ${}_{p+1}
  F_{p}$, but if the $\lambda_{n}^{(\alpha)}$ have a complicated
  structure or if only their numerical values are available, numerical
  summation techniques have to be employed. It is shown that certain
  nonlinear sequence transformations -- in particular the so-called delta
  transformation [E.\ J.\ Weniger, Comput.\ Phys.\ Rep.\ \textbf{10}, 189
  -- 371 (1989), Eq.\ (8.4-4)] -- are able to sum the divergent series
  occurring in this context effectively. As a physical application of the
  results of this article, the legitimacy of the rearrangement of certain
  one-range addition theorems for Slater-type functions [I.\ I.\
  Guseinov, Phy. Rev. A \textbf{22}, 369 - 371 (1980); Int. J. Quantum
  Chem. \textbf{81}, 126 - 129 (2001); Int.  J. Quantum Chem.
  \textbf{90}, 114 - 118 (2002)] is investigated.
\end{abstract}

\begin{quote}
  {\bfseries PACS numbers:} 02.30.Gp, 02.30.Lt, 02.30.Mv, 02.60.-x

  {\bfseries AMS classification scheme numbers:}  30B10, 30B40, 33C45,
  40A05, 65B10

  \noindent {\bfseries Keywords:} Generalized Laguerre polynomials,
  Laguerre series, analyticity, divergent series, summability, nonlinear
  sequence transformations
\end{quote}

\newpage

\tableofcontents

\newpage


\typeout{==> Section: Introduction}
\section{Introduction}
\label{Sec:Introduction}

The generalized Laguerre polynomials $L_{n}^{(\alpha)} (z)$, whose most
relevant properties are reviewed in Section
\ref{Sec:GeneralizedLaguerrePolynomials}, are a very important class of
orthogonal polynomials with numerous mathematical and physical
applications. There is an extensive literature both on their mathematical
properties as well as on their applications, and any attempt of providing
a reasonably complete bibliography would be hopeless. Let me just mention
that the radial parts of bound state hydrogen eigenfunctions and of
several other physically relevant complete and orthonormal sets of
function $f \colon \mathbb{R}^{3} \to \mathbb{C}$ are essentially
generalized Laguerre polynomials (see for example \cite[Sections IV and
V]{Weniger/1985} and references therein).

As discussed in more detail in Section
\ref{Sec:GeneralizedLaguerrePolynomials}, the generalized Laguerre
polynomials $L_{n}^{(\alpha)} (z)$ form a complete orthogonal polynomial
system in the weighted Hilbert space $L^{2}_{z^{\alpha} \mathrm{e}^{-z}}
\bigl( [0, \infty)\bigr)$ defined by (\ref{HilbertL^2_Lag}), which is
based on an integration over the interval $[0, \infty)$ involving the
weight function $w (z) = z^{\alpha} \exp (-z)$. Accordingly, the topic of
this article are infinite Laguerre series of the following type:
\begin{subequations}
  \label{f_Exp_GLag}
  \begin{align}
    \label{f_Exp_GLag_a}
    f (z) & \; = \; \sum_{n=0}^{\infty} \,
    \lambda_{n}^{(\alpha)} \, L_{n}^{(\alpha)} (z) \, ,
    \\
    \label{f_Exp_GLag_b}
    \lambda_{n}^{(\alpha)} & \; = \; \frac{n!}{\Gamma (\alpha+n+1)} \,
    \int_{0}^{\infty} \, z^{\alpha} \, \mathrm{e}^{-z} \,
    L_{n}^{(\alpha)} (z) \, f (z) \, \mathrm{d} z \, .
  \end{align}
\end{subequations}
The expansion coefficients $\lambda_{n}^{(\alpha)}$ are essentially inner
products, utilizing the orthogonality and the completeness of the
generalized Laguerre polynomials in $L^{2}_{z^{\alpha} \mathrm{e}^{-z}}
\bigl( [0, \infty)\bigr)$.

Not all expansions in terms of generalized Laguerre polynomials can be
interpreted in a Hilbert space setting. There are expansions of special
functions in terms of generalized Laguerre polynomials with variable and
index-dependent superscripts. A simple example is the generating function
\cite[p.\ 242]{Magnus/Oberhettinger/Soni/1966}
\begin{equation}
  \label{GLag_Gen_5}
  \sum_{n=0}^{\infty} \, L_{n}^{(\alpha-n)} (x) \, t^{n} \; = \;
    \mathrm{e}^{-xt} \, (1+t)^{\alpha} \, , \qquad \vert t \vert < 1 \, .
\end{equation}
Expansions of that kind cannot be derived via a straightforward
application of the orthogonality of the generalized Laguerre polynomials.
In the article by S\'{a}nchez-Ruiz, L\'{o}pez-Art\'{e}z, and Dehesa
\cite{Sanchez-Ruiz/Lopez-Artez/Dehesa/2003}, transformation formulas for
generalized hypergeometric series were employed for the construction of
Laguerre expansions with variable superscripts. However, nonorthogonal
Laguerre expansions of that kind are not the topic of this article.

It is generally accepted that orthogonal expansions are extremely useful
mathematical tools and that they have many highly advantageous features.
This is, however, not the whole truth, in particular if we want to
approximate functions.  There are situations, in which alternative
representations are (much) more useful. Hilbert space theory only
guarantees that an orthogonal expansion converges in the mean with
respect to the corresponding norm, but not necessarily pointwise or even
uniformly.  Therefore, orthogonal expansions are not necessarily the best
choice if we are interested in \emph{local} properties of functions.

Convergence in the mean is -- although completely satisfactory for many
purposes such as the evaluation of integrals -- a comparatively weak form
of convergence. In practice, it is therefore often desirable or even
necessary to construct alternative representations possessing more
convenient properties (see also the discussion in
\cite{Ford/Pennline/2007}).

A very desirable feature of functions $f \colon \mathbb{C} \to
\mathbb{C}$ is \emph{analyticity} in the sense of complex analysis. This
means that a function $f$ can be represented in a neighborhood of the
origin by a convergent power series,
\begin{equation}
  \label{PowSer_f}
  f (z) \; = \; \sum_{n=0}^{\infty} \, \gamma_{n} \, z^{n} \, ,
\end{equation}
and that the coefficients $\gamma_{n}$ of this power series essentially
correspond to the derivatives of $f$ at $z=0$.

If the function $f$ defined by the Laguerre series (\ref{f_Exp_GLag}) is
explicitly known, it is usually not too difficult to decide whether $f$
is analytic or not, and if $f$ is analytic it is normally not too
difficult to construct at least the leading terms of its power series
expansion (\ref{PowSer_f}).

Unfortunately, the situation is not always so good, and it can even
happen that only the numerical values of a finite number of expansion
coefficients $\lambda_{n}^{(\alpha)}$ in (\ref{f_Exp_GLag}) are known. In
such a case, it would certainly be helpful if we could relate the
properties of the expansion coefficients $\lambda_{n}^{(\alpha)}$ -- in
particular their decay rate and their sign pattern -- to the analyticity
of the function $f (z)$ defined by the Laguerre series
(\ref{f_Exp_GLag}). It would also be helpful if we could construct at
least the leading power series coefficients $\gamma_{n}$ -- either
exactly or in an approximate sense -- from the coefficients
$\lambda_{n}^{(\alpha)}$ of the Laguerre series (\ref{f_Exp_GLag}). .

Pollard \cite{Pollard/1947a}, Sz\'{a}sz and Yeardley
\cite{Szasz/Yeardley/1958}, and Rusev \cite{Rusev/1983} had investigated
Laguerre expansions of the type of (\ref{f_Exp_GLag}) of analytic
functions and analyzed their regions of convergence. However, the inverse
problem -- the formulation of criteria which guarantee that a Laguerre
expansion of the type of (\ref{f_Exp_GLag}) represents a function
analytic in a neighborhood of the origin -- seems to be essentially
unexplored. I am only aware of short remarks by Gottlieb and Orszag
\cite[p.\ 42]{Gottlieb/Orszag/1977} and by Doha \cite[p.\
5452]{Doha/2003}, respectively, who stated without detailed proof that
such a Laguerre series converges faster than algebraically if the
function under consideration is analytic at the origin.

Many Laguerre series are known which seem to confirm the claim of
Gottlieb and Orszag and of Doha, respectively. The probably most simple
example is the well known generating function $(1-t)^{-\alpha-1} \exp
\bigl( zt/[t-1] \bigr)$ of the Laguerre polynomials. For $\vert t \vert <
1$, it is an \emph{entire} function, and the coefficients of its Laguerre
series (\ref{GLag_Gen_1}) decay exponentially for $\vert t \vert < 1$.

Another example is the Laguerre series (\ref{GenPow2GLag}) for the
general power function $z^{\rho}$ with nonintegral $\rho \in \mathbb{R}
\setminus \mathbb{N}_0$. Obviously, $z^{\rho}$ with $\rho \in \mathbb{R}
\setminus \mathbb{N}_0$ is not analytic at $z=0$. The series coefficients
in (\ref{GenPow2GLag}) decay algebraically and possess for sufficiently
large indices the same sign. As shown in Section
\ref{Sec:LagSerGeneralPowerFunction}, it is nevertheless possible to
transform the Laguerre series (\ref{GenPow2GLag}) for $z^{\rho}$ with
$\rho \in \mathbb{R} \setminus \mathbb{N}_0$ to the \emph{formal} power
series (\ref{Chk_x^m_GlagPol_6}). However, the formal power series
(\ref{Chk_x^m_GlagPol_6}) is not a mathematically meaningful object: For
sufficiently large indices, the power series coefficients in
(\ref{Chk_x^m_GlagPol_6}) are according to
(\ref{Lim_Ser_Chk_x^m_GlagPol_6}) all infinite in magnitude. Thus, the
Laguerre series (\ref{GenPow2GLag}) for $z^{\rho}$ with $\rho \in
\mathbb{R} \setminus \mathbb{N}_0$ also agrees with the claim of Gottlieb
and Orszag and of Doha, respectively.

If, however, we augment the Laguerre series (\ref{GenPow2GLag}) for
$z^{\rho}$ by an additional alternating sign $(-1)^{n}$, we obtain the
Laguerre series (\ref{LagSer_AltSerCoeffs}). Its transformation to a
power series yields according to (\ref{Rearr_G_Exp_GLag_4}) a confluent
hypergeometric function ${}_{1} F_{1} \bigl( -\rho; \alpha+1; z/2
\bigr)$, which for $\Re (\alpha) > -1$ is an \emph{entire} function. The
coefficients in the Laguerre series (\ref{LagSer_AltSerCoeffs}) decay
algebraically in magnitude. This implies that the transformation of
(\ref{LagSer_AltSerCoeffs}) to a power series leads to power series
coefficients $\gamma_{n}$ that are represented by divergent series
expansions. However, the alternating sign $(-1)^{n}$ in
(\ref{LagSer_AltSerCoeffs}) makes it possible to associate finite values
to these divergent series, i.e., the series expansions for the power
series coefficients are now \emph{summable}. Thus, we obtain a function,
which is analytic at $z=0$ and which does not agree with the claim of
Gottlieb and Orszag and of Doha, respectively.

The example of the closely related Laguerre series (\ref{GenPow2GLag})
and (\ref{LagSer_AltSerCoeffs}) shows that it makes a huge difference if
Laguerre series coefficients $\lambda_{n}^{(\alpha)}$, that decay
algebraically in magnitude, ultimately have strictly alternating or
strictly monotone signs. Thus, the claim of Gottlieb and Orszag and of
Doha, respectively, is imprecise and ignores the pivotal role of
divergent, but summable series expansions for the power series
coefficients $\gamma_{n}$ in (\ref{PowSer_f}). Some basic facts about the
summation of divergent series are reviewed in Appendix
\ref{App:DivergentSeries}.

A detailed investigation of the analyticity of functions represented by
Laguerre expansions of the type of (\ref{f_Exp_GLag}) is the topic of
this article. The central result is the transformation formula
(\ref{Rearr_f_Exp_GLag}), which yields a \emph{formal} power series
expansion for a function represented by a Laguerre expansion. The crucial
question is whether the inner $\mu$ series in (\ref{Rearr_f_Exp_GLag})
for the power series coefficients $\gamma_{n}$ converge. If these $\mu$
series do not converge and are also not summable to something finite,
then the function under consideration is not analytic at the origin and
the resulting formal power series that is mathematically meaningless
since contains coefficients $\gamma_{n}$ which are infinite in magnitude.
Consequently, it is comparatively easy to relate the analyticity of a
function $f$ represented by a Laguerre series at the origin $z=0$ with
the decay rate and the sign pattern of its coefficients
$\lambda_{n}^{(\alpha)}$.

Analytical manipulations can only lead to closed form expressions for the
power series coefficients $\gamma_{n}$ if the Laguerre series
coefficients $\lambda_{n}^{(\alpha)}$ possess a sufficiently simple
structure. This undeniable fact may create the false impression that the
formalism developed in this article is restricted to Laguerre series with
very simple coefficients $\lambda_{n}^{(\alpha)}$, and that this
formalism is at best suited for an alternative rederivation of known
generating functions of the generalized Laguerre polynomials.

As discussed in Section \ref{Sec:ComputationalApproaches}, it is,
however, often possible to construct \emph{numerical approximations} to
the leading power series coefficients $\gamma_{n}$ from the coefficients
$\lambda_{n}^{(\alpha)}$ of the Laguerre series, even if the inner $\mu$
series in (\ref{Rearr_f_Exp_GLag}) diverges. The necessary summations can
be done effectively with the help of certain nonlinear sequence
transformations -- Wynn's epsilon algorithm (\ref{eps_al}) and the two
Levin-type transformations (\ref{dLevTr}) and (\ref{dWenTr}) -- whose
basic properties are reviewed in Appendix
\ref{App:SequenceTransformations}. While Wynn's epsilon algorithm
(\ref{eps_al}), which produces Pad\'{e} approximants if the input data
are the partial sums of a power series, is now fairly well known among
(applied) mathematicians and theoretical physicists, the Levin-type
transformations (\ref{dLevTr}) and (\ref{dWenTr}), which in some cases
were found to be remarkably powerful, are not nearly as well known as
they deserve to be.

As a physical application of the mathematical formalism developed in this
article, the legitimacy of rearrangements of certain one-range addition
theorems for Slater-type functions with in general nonintegral principal
quantum numbers performed by Guseinov
\cite{Guseinov/1980a,Guseinov/2001a,Guseinov/2002c} is investigated in
Section \ref{Sec:GuseinovsRearrangedOne-RangeAdditionTheorems}. It is
shown that the one-center limits of Guseinov's rearranged addition
theorems do not exist if the principal quantum numbers of the Slater-type
functions are nonintegral.

\typeout{==> Section: Generalized Laguerre Polynomials}
\section{Generalized Laguerre Polynomials}
\label{Sec:GeneralizedLaguerrePolynomials}

The generalized Laguerre polynomials $L_{n}^{(\alpha)} (z)$ with $n \in
\mathbb{N}_{0}$ and $\Re (\alpha) > - 1$ are orthogonal polynomials
associated with the integration interval $[0, \infty)$ and the weight
function $w (z) = z^{\alpha} \exp (-z)$.

In this article, the mathematical notation for generalized Laguerre
polynomials (see for example \cite[Chapter
5.5]{Magnus/Oberhettinger/Soni/1966}) is used. Additional conventions, as
used predominantly in older books and articles on quantum theory, were
discussed by Kaijser and Smith \cite[Footnote 1 on p.\
48]{Kaijser/Smith/1977} and also in \cite[Section 4]{Weniger/2007b}.

The generalized Laguerre polynomials possess an
explicit expression as a terminating confluent hypergeometric series
\cite[p.\ 240]{Magnus/Oberhettinger/Soni/1966},
\begin{equation}
  \label{GLag_1F1}
L_{n}^{(\alpha)} (z) \; = \;
\frac{(\alpha+1)_n}{n!} \, {}_1 F_1 (-n; \alpha+1; z) \, ,
\end{equation}
and they can also be defined via their Rodrigues relationship \cite[p.\
241]{Magnus/Oberhettinger/Soni/1966}
\begin{equation}
  \label{GLag_Rodrigues}
L_{n}^{(\alpha)} (z) \; = \; z^{-\alpha} \, \frac{\mathrm{e}^{z}}{n!} \,
\frac{\mathrm{d}^n}{\mathrm{d} z^n} \,
\bigl[ \mathrm{e}^{-z} z^{n+\alpha} \bigr] \, .
\end{equation}
In the special case $\alpha=0$, it is common to drop the superscript
since we obtain the usual Laguerre polynomials \cite[p.\
239]{Magnus/Oberhettinger/Soni/1966}:
\begin{equation}
  L_n (z) \; = \; L_{n}^{(0)} (z) \, , \qquad n \in \mathbb{N}_0 \, .
\end{equation}

Either the explicit expression (\ref{GLag_1F1}) or the Rodrigues
relationship (\ref{GLag_Rodrigues}) can be used to define generalized
Laguerre polynomials with essentially arbitrary complex values of the
superscript $\alpha$. But in the orthogonality relationship of the
generalized Laguerre polynomials \cite[p.\
241]{Magnus/Oberhettinger/Soni/1966},
\begin{equation}
  \label{GLag_Orthogonality}
\int_{0}^{\infty} \, z^{\alpha} \, \mathrm{e}^{-z} \,
L_{m}^{(\alpha)} (z) \, L_{n}^{(\alpha)} (z) \, \mathrm{d} z \; = \;
\frac{\Gamma (\alpha+n+1)}{n!} \, \delta_{m n} \, ,
\qquad m, n \in \mathbb{N}_{0 \, ,}
\end{equation}
it is necessary to assume $\Re (\alpha) > - 1$ because otherwise this
integral does not exist. However, as discussed by Kochneff
\cite{Kochneff/1995}, this condition can be relaxed if this and related
integrals are reinterpreted as Hadamard finite part integrals.

In the vast majority of calculations involving generalized Laguerre
polynomials, the superscript $\alpha$ is real and also positive. In
practice, more general complex values of $\alpha$ play a negligible rule.
Therefore, in the following text it is always tacitly assumed that all
integrals exist in the ordinary sense, which avoids the complications of
Hadamard finite part integrals, and that the the superscript $\alpha$ of
a generalized Laguerre polynomial is a real number satisfying $\alpha >
-1$.

The orthogonality relationship (\ref{GLag_Orthogonality}) suggests the
introduction of the following inner product containing the weight
function $w (z) = z^{\alpha} \exp (-z)$ for functions $f, g \colon
\mathbb{C} \to \mathbb{C}$:
\begin{equation}
  \label{Def:InnerProduct_Lag}
  ( f \vert g )_{z^{\alpha} \mathrm{e}^{-z}, 2} \; = \;
  \int_{0}^{\infty} \, \, z^{\alpha} \mathrm{e}^{-z} \, [f (z)]^{*} \,
  g (z) \, \mathrm{d} z \, .
\end{equation}
This inner product gives rise to the norm
\begin{equation}
  \label{Def:Norm_Lag}
  \Vert f \Vert_{z^{\alpha} \mathrm{e}^{-z}, 2} \; = \;
  \sqrt{( f \vert f )_{z^{\alpha} \mathrm{e}^{-z}, 2}} \, .
\end{equation}
We then obtain the following Hilbert space of square integrable
functions:
\begin{align}
  \label{HilbertL^2_Lag}
  L^{2}_{z^{\alpha} \mathrm{e}^{-z}} \bigl( [0, \infty)\bigr) & \; = \;
  \Bigl\{ f \colon \mathbb{C} \to \mathbb{C} \Bigm\vert \,
  \int_{0}^{\infty} \, z^{\alpha} \, \mathrm{e}^{-z} \, \vert f (z)
  \vert^2 \, \mathrm{d} z < \infty \Bigr\}
  \notag \\[1\jot]
  & \; = \; \bigl\{ f \colon \mathbb{C} \to \mathbb{C} \big\vert \, \Vert
  f \Vert_{z^{\alpha} \mathrm{e}^{-z}, 2} < \infty \bigr\} \, .
\end{align}

The completeness of the generalized Laguerre polynomials in this weighted
Hilbert space is a classic result of mathematical analysis (see for
example \cite[p.\ 33]{Higgins/1977}, \cite[pp.\ 349 - 351]{Sansone/1977},
or \cite[pp.\ 235 - 238]{Tricomi/1970}). Thus, the normalized polynomials
\begin{equation}
  \label{Def:GenLag_Normalized}
  \mathcal{L}_{n}^{(\alpha)} (z) \; = \;
  \left[ \frac{n!}{\Gamma (\alpha+n+1)} \right]^{1/2} \,
  L_{n}^{(\alpha)} (z) \, ,
  \qquad n \in \mathbb{N}_0 \, ,
\end{equation}
are a complete and orthonormal polynomial system in $L^{2}_{z^{\alpha}
  \mathrm{e}^{-z}} \bigl( [0, \infty)\bigr)$, satisfying
\begin{equation}
  \label{Orthonorm_CalL}
  \int_{0}^{\infty} \, z^{\alpha} \, \mathrm{e}^{-z} \,
  \mathcal{L}_{m}^{(\alpha)} (z) \, \mathcal{L}_{n}^{(\alpha)} (z) \,
  \mathrm{d} z \; = \; \delta_{m n} \, .
\end{equation}
Accordingly, every $f \in L^{2}_{z^{\alpha} \mathrm{e}^{-z}} \bigl( [0,
\infty)\bigr)$ can be expanded in terms of the normalized polynomials
$\bigl\{ \mathcal{L}_{n}^{(\alpha)} (z) \bigr\}_{n=0}^{\infty}$:
\begin{subequations}
  \label{f_Exp_GenLag_Normalized}
  \begin{align}
    \label{f_Exp_GenLag_Normalized_a}
    f (z) & \; = \; \sum_{n=0}^{\infty} \, \mathcal{C}_{n}^{(\alpha)} \,
    \mathcal{L}_{n}^{(\alpha)} (z) \, ,
    \\
    \label{f_Exp_GenLag_Normalized_b}
    \mathcal{C}_{n}^{(\alpha)} & \; = \; \int_{0}^{\infty} \,
    z^{\alpha} \, \mathrm{e}^{-z} \,
    \mathcal{L}_{n}^{(\alpha)} (z) \, f (z) \, \mathrm{d} z \, .
  \end{align}
\end{subequations}
This expansion in terms of the normaized polynomials
$\mathcal{L}_{n}^{(\alpha)} (z)$, which converges in the mean with
respect to the norm (\ref{Def:Norm_Lag}), is nothing but the Laguerre
series (\ref{f_Exp_GLag}) in disguise. In agreement with
(\ref{Def:GenLag_Normalized}), we only have to choose
\begin{equation}
  \label{CalC_2_lambda}
  \mathcal{C}_{n}^{(\alpha)} \; = \;
  \left[ \frac{\Gamma (\alpha+n+1)}{n!} \right]^{1/2} \,
  \lambda_{n}^{(\alpha)}
\end{equation}
to see that the expansions (\ref{f_Exp_GLag}) and
(\ref{f_Exp_GenLag_Normalized}) are identical.

Let me emphasize once more that Laguerre expansions converge in general
only in the mean, but not necessarily pointwise (see for example
\cite{Askey/Wainger/1965}). Additional conditions, which a function has
to satisfy in order to guarantee that its Laguerre expansion also
converges pointwise, were discussed by Szeg\"{o} \cite[Theorem 9.1.5 on
p.\ 246]{Szegoe/1967}.

Hilbert space theory can be used to derive sufficient criteria, which the
coefficients $\lambda_{n}^{{\alpha}}$ or $\mathcal{C}_{n}^{(\alpha)}$
have to satisfy in order to guarantee that the equivalent expansions
(\ref{f_Exp_GLag}) and (\ref{f_Exp_GenLag_Normalized}) converge in the
mean. The basic requirement is that the norm (\ref{Def:Norm_Lag}) of both
$f$ and its Laguerre expansion (\ref{f_Exp_GenLag_Normalized}) must be
finite:
\begin{equation}
  \label{ConvCond_f_Exp_GenLag_Normalized}
  \Vert f \Vert_{z^{\alpha} \mathrm{e}^{-z}, 2} \; = \;
  \biggl\{ \sum_{n=0}^{\infty} \,
  \bigl\vert \mathcal{C}_{n}^{(\alpha)} \bigr\vert^{2} \biggr\}^{1/2}
  \; < \; \infty \, .
\end{equation}
It follows from (\ref{CalC_2_lambda}) that this condition can also be
reformulated as follows:
\begin{equation}
  \label{MeanConvCond_f}
  \Vert f \Vert_{z^{\alpha} \mathrm{e}^{-z}, 2}
  \; = \; \biggl\{ \sum_{n=0}^{\infty} \,
  \frac{\Gamma (\alpha+n+1)}{n!} \,
  \bigl\vert \lambda_{n}^{(\alpha)} \bigr\vert^{2} \biggr\}^{1/2}
  \; < \; \infty \, .
\end{equation}

A sufficient condition, which guarantees that the infinite series
$\sum_{n=0}^{\infty} \bigl\vert \mathcal{C}_{n}^{(\alpha)}
\bigr\vert^{2}$ in (\ref{ConvCond_f_Exp_GenLag_Normalized}) converges, is
\begin{equation}
  \label{BigO_C}
  \bigl\vert \mathcal{C}_{n}^{(\alpha)} \bigr\vert^{2}
  \; \sim \; n^{-1-\epsilon} \, ,
  \qquad n \to \infty \, , \quad \epsilon > 0 \, .
\end{equation}
It follows from (\ref{CalC_2_lambda}) and (\ref{BigO_C}) that the
Laguerre series (\ref{f_Exp_GLag}) converges in the mean if its
coefficients $\lambda_{n}^{(\alpha)}$ satisfy
\begin{equation}
  \frac{\Gamma (\alpha+n+1)}{n!} \,
  \bigl\vert \lambda_{n}^{(\alpha)} \bigr\vert^{2}
  \; \sim \; n^{-1-\epsilon} \, ,
  \qquad n \to \infty \, , \quad \epsilon > 0 \, .
\end{equation}
Now, we need the asymptotic approximation \cite[Eq.\ (6.1.47) on p.\
257]{Abramowitz/Stegun/1972}
\begin{equation}
  \label{AsyGammaRatio}
  \Gamma(z+a)/\Gamma(z+b) \; = \; z^{a-b} \,
  \bigl[ 1 + \mathrm{O} (1/z) \bigr] \, , \qquad z \to \infty \, ,
\end{equation}
which is the leading term of an asymptotic expansion first derived by
Tricomi and Erd\'{e}lyi \cite{Tricomi/Erdelyi/1951} that holds -- as
emphasized by Olver \cite[p.\ 119]{Olver/1997a} -- without restrictions
on $a, b \in \mathbb{C}$. Thus, we obtain:
\begin{equation}
  \label{Asy_Coeff_LagSer}
  \frac{\Gamma (\alpha+n+1)}{n!} \; \sim \; n^{\alpha}
  \, + \,  \mathrm{O} \bigl( n^{\alpha-1} \bigr) \, ,
  \qquad n \to \infty \, .
\end{equation}
This translates to the sufficient convergence condition that $\bigl\vert
\lambda_{n}^{(\alpha)} \bigr\vert$ must decay like
\begin{equation}
  \label{ConvCon_LambdaSer}
  \bigl\vert \lambda_{n}^{(\alpha)} \bigr\vert \; \sim \;
  n^{-[\alpha+\epsilon+1]/2} \, ,
  \qquad n \to \infty \, , \qquad \epsilon > 0 \, ,
\end{equation}
or faster. If this asymptotic condition is satisfied, then the Laguerre
series (\ref{f_Exp_GLag}) for $f (z)$ converges in the mean with respect
to the norm (\ref{Def:Norm_Lag}) of the Hilbert space
$L^{2}_{\mathrm{e}^{-z} z^{\alpha}} \bigl( [0, \infty) \bigr)$.

The sufficient convergence condition (\ref{ConvCon_LambdaSer}) shows that
the coefficients $\lambda_{n}^{(\alpha)}$ can decay extremely slowly,
which translates to a possibly extremely slow convergence of Laguerre
series of the type of (\ref{f_Exp_GLag}). Thus, it would be overly
optimistic to assume that all Laguerre series are necessarily
computationally useful.

\typeout{==> Section: The Laguerre Series for the General Power Function}
\section{The Laguerre Series for the General Power Function}
\label{Sec:LagSerGeneralPowerFunction}

In the theory of orthogonal expansions, which converge in the mean with
respect to the norm $\Vert \cdot \Vert$ of some Hilbert space
$\mathcal{H}$, the decisive criterion is that the function, which is to
be expanded, has to belong to $\mathcal{H}$. Consequently, any function
$f \colon \mathbb{C} \to \mathbb{C}$ with finite norm
(\ref{Def:Norm_Lag}) belongs to the weighted Hilbert space
$L^{2}_{z^{\alpha} \mathrm{e}^{-z}} \bigl( [0, \infty)\bigr)$ and can be
expanded in terms of generalized Laguerre polynomials according to
(\ref{f_Exp_GLag}). The resulting expansion (\ref{f_Exp_GLag}) converges
in the mean with respect to the norm (\ref{Def:Norm_Lag}) of the Hilbert
space $L^{2}_{z^{\alpha} \mathrm{e}^{-z}} \bigl( [0, \infty)\bigr)$.

The existence and convergence of a Laguerre series of the type of
(\ref{f_Exp_GLag}) does not allow any conclusion about the analyticity of
the corresponding function. The weighted Hilbert space $L^{2}_{z^{\alpha}
  \mathrm{e}^{-z}} \bigl( [0, \infty)\bigr)$ contains many functions that
are obviously not analytic at the origin. Therefore, attempts of
constructing a power series from the Laguerre series can easily lead to a
disaster. The complications, which can occur in this context, can be
demonstrated via the general power function $z^{\rho}$ with nonintegral
$\rho \in \mathbb{R} \setminus \mathbb{N}_{0}$.

In integrals over the positive real semi-axis, the weight function
$z^{\alpha} \mathrm{e}^{-z}$ becomes at least for $\Re (\alpha) > 0$ very
small both as $z \to 0$ and as $z \to \infty$. Consequently, $z^{\alpha}
\mathrm{e}^{-z}$ suppresses in the inner product
(\ref{Def:InnerProduct_Lag}) the contribution of the remaining integrand
for small and large arguments. Therefore, the general power function
$z^{\rho}$ with $\rho \in \mathbb{R} \setminus \mathbb{N}_{0}$ possesses
at least for sufficiently large values of $\alpha$ an expansion of the
type of (\ref{f_Exp_GLag}) in terms of generalized Laguerre polynomials.
The existence of this expansion is guaranteed if $z^{\rho}$ belongs to
the Hilbert space $L^{2}_{z^{\alpha} \mathrm{e}^{-z}} \bigl( [0,
\infty)\bigr)$, or equivalently, if
\begin{equation}
  \left\vert \int_{0}^{\infty} \, z^{\alpha+2\rho} \, \mathrm{e}^{-z} \,
    \mathrm{d} z \right\vert \; < \; \infty
\end{equation}
holds. Thus, we have to require that $\alpha+2\rho > - 1$ holds.

For the construction of a Laguerre series for $z^{\rho}$, we only need the
integral \cite[Eq.\ (7.414.7) on p.\ 850]{Gradshteyn/Rhyzhik/1994}
\begin{align}
  \label{Gr_7.414.7}
  & \int_{0}^{\infty} \, \mathrm{e}^{-st} \, t^{\beta} \,
  L_{n}^{(\alpha)} (t) \, \mathrm{d} t
  \notag \\
  & \qquad \; = \; \frac{\Gamma (\beta+1) \, \Gamma (\alpha+n+1)}{n! \,
    \Gamma (\alpha+1)} \, s^{-\beta-1} \, {}_2 F_1 (-n, \beta+1;
  \alpha+1; 1/s) \, ,
  \notag \\
  & \qquad \qquad \Re (\beta) > - 1 \, , \qquad \Re (s) > 0 \, ,
\end{align}
to obtain after some essentially straightforward algebra \cite[Eq.\ (16)
on p.\ 214]{Erdelyi/Magnus/Oberhettinger/Tricomi/1953b}:
\begin{equation}
  \label{GenPow2GLag}
  z^{\rho} \; = \; \frac{\Gamma (\rho+\alpha+1)}{\Gamma (\alpha+1)} \,
  \sum_{n=0}^{\infty} \, \frac{(-\rho)_n}{(\alpha+1)_n} \,
  L_{n}^{(\alpha)} (z) \, ,
  \qquad \rho \in \mathbb{R} \setminus \mathbb{N}_0 \, ,
  \qquad \alpha+2\rho > - 1 \, .
\end{equation}
If we set $\rho=m$ with $m \in \mathbb{N}_0$, then the infinite series on
the right-hand side terminates because of the Pochhammer symbol $(-m)_n$
and we obtain the following finite sum (see for example \cite[Eq.\ (2) on
p.\ 207]{Rainville/1971}):
\begin{equation}
  \label{IntPow_GlagPol}
  z^{m} \; = \;
  (\alpha+1)_m \, \sum_{n=0}^{m} \, \frac{(-m)_n}{(\alpha+1)_n}
  \, L_{n}^{(\alpha)} (z) \, ,
  \qquad m \in \mathbb{N}_0 \, , \qquad \alpha+2m > - 1 \, .
\end{equation}

Although intimately related, there are nevertheless some fundamental
differences between the two Laguerre expansions (\ref{GenPow2GLag}) and
(\ref{IntPow_GlagPol}). The finite sum formula (\ref{IntPow_GlagPol}) is
a relationship among polynomials and therefore certainly valid pointwise
for arbitrary $z \in \mathbb{C}$ as well as analytic at the origin in the
sense of complex analysis.

In the case of the infinite series expansion (\ref{GenPow2GLag}), we only
know that it converges in the mean with respect to the norm
(\ref{Def:Norm_Lag}), but we have no \emph{a priori} reason to assume
that this expansion might converge pointwise for arbitrary $\rho \in
\mathbb{R} \setminus \mathbb{N}_{0}$. Moreover, $z^{\rho}$ with $\rho \in
\mathbb{R} \setminus \mathbb{N}_{0}$ is not analytic at the origin.

The validity of the finite sum (\ref{IntPow_GlagPol}) can be checked by
inserting the explicit expression (\ref{GLag_1F1}) for the generalized
Laguerre polynomials. If we then rearrange the order of the two nested
sums, we obtain after some algebra:
\begin{equation}
  \label{Chk_IntPow_GlagPol_6}
  z^{m} \; = \; (\alpha+1)_m \, \sum_{k=0}^{m} \, (-1)^k \,
  \frac{(-m)_k}{(\alpha+1)_k} \, \frac{z^k}{k!} \, \sum_{\nu=0}^{m-k} \,
  (-1)^{\nu} \, \binom{m-k}{\nu} \, .
\end{equation}
Next, we use the relationship \cite[Eq.\ (3.1.7) on p.\
10]{Abramowitz/Stegun/1972}:
\begin{equation}
  \sum_{k=0}^{n} \, (-1)^k \, \binom{n}{k} \; = \; \delta_{n 0} \, ,
  \qquad n \in \mathbb{N}_0 \, ,
\end{equation}
for binomial sums, which shows that the inner sum in
(\ref{Chk_IntPow_GlagPol_6}) vanishes unless we have $k=m$. Thus, we only
need $(-1)^m (-m)_m = (1)_m = m!$ to arrive at the trivial identity
$z^{m}=z^{m}$ which proves the correctness of (\ref{IntPow_GlagPol}).

In the case of the infinite series (\ref{GenPow2GLag}) for $z^{\rho}$ with
$\rho \in \mathbb{R} \setminus \mathbb{N}_{0}$, we can also insert the
explicit expression (\ref{GLag_1F1}) for the generalized Laguerre
polynomials into it and rearrange the order of summations. We then obtain
after some algebra:
\begin{equation}
  \label{Chk_x^m_GlagPol_6}
  z^{\rho} \; = \; \frac{\Gamma (\rho+\alpha+1)}{\Gamma (\alpha+1)}
  \, \sum_{k=0}^{\infty} \, (-1)^k \, \frac{(-\rho)_k}{(\alpha+1)_k} \,
  \frac{z^k}{k!} \, {}_1 F_0 (k-\rho; 1) \, .
\end{equation}
Superficially, it looks as if we succeeded in constructing a power series
expansion for the in general nonintegral power $z^{\rho}$. However, the
generalized hypergeometric series ${}_1 F_0$ with unit argument is the
limiting case $z \to 1$ of the so-called binomial series \cite[p.\
38]{Magnus/Oberhettinger/Soni/1966}:
\begin{equation}
  \label{BinomSer}
  {}_1 F_0 (a; z) \; = \;
  \sum_{m=0}^{\infty} \, \binom{-a}{m} (-z)^m \; = \; (1-z)^{-a} \, ,
  \qquad \vert z \vert < 1 \, .
\end{equation}
If we set $a=k-\rho$ with $k \in \mathbb{N}_{0}$ and $\rho \in \mathbb{R}
\setminus \mathbb{N}_{0}$, we obtain for the hypergeometric series ${}_1
F_0$ in (\ref{Chk_x^m_GlagPol_6}):
\begin{equation}
  \label{Lim_Ser_Chk_x^m_GlagPol_6}
  {}_1 F_0 (k-\rho; 1) \; = \; \lim_{z \to 1} \, (1-z)^{\rho-k} \; = \;
  \begin{cases}
    \infty \, , \qquad \rho < 0 \, , \\
    0 \, , \qquad \; \, k < \rho \ge 0 \, , \\
    \infty \, , \qquad k > \rho \ge 0 \, .
  \end{cases}
\end{equation}
Thus, disaster struck and the power series (\ref{Chk_x^m_GlagPol_6})
contains infinitely many series coefficients that are infinite in
magnitude.

I am aware of an article by Villani \cite{Villani/1972} who tried to make
sense of perturbation expansions with divergent terms. Since, however,
this is the only article on this topic, which I am aware of, I am tempted
to believe that Villani's approach was not overly fertile.  Therefore, I
will stick to the usual mathematical convention that a power series with
coefficients, that are infinite in magnitude, does not exist as a
mathematically meaningful object.

It is important to note that series with \emph{divergent terms} and
\emph{divergent series} are not the same. In the case of divergent
series, all terms are finite, but the conventional process of adding up
the terms successively does not lead to a convergent result.
Nevertheless, it is often possible to associate a finite value to a
divergent series with the help of a suitable summation technique. As
reviewed in Appendix \ref{App:DivergentSeries}, divergent series and
their summation have been and to some extent still are a highly
controversial topic. The summation of divergent series plays a major role
in Sections (\ref{Sec:AlgebraicallyDecayingSeriesCoefficients}) and
(\ref{Sec:ComputationalApproaches}).

Since the hypergeometric series ${}_{1} F_{0} (k-\rho; 1)$ in
(\ref{Chk_x^m_GlagPol_6}) does not converge for all indices $k$, the
Laguerre series (\ref{GenPow2GLag}) for $z^{\rho}$ with $\rho \in
\mathbb{R} \setminus \mathbb{N}_0$ cannot be reformulated as a power
series in $z$ by interchanging the order of summations. Of course, this
is a mathematical necessity: The general power function $z^{\rho}$ with
$\rho \in \mathbb{R} \setminus \mathbb{N}_0$ is not analytic at the
origin, which implies that a power series about $z=0$ cannot exist.

So far, the analysis of this Section has only produced obvious results
and no new insight: The integral power $z^{m}$ with $m \in
\mathbb{N}_{0}$ is analytic at the origin. Consequently, it must be
possible to reformulate its \emph{finite} Laguerre expansion
(\ref{IntPow_GlagPol}) as a polynomial in $z$. In contrast, the
nonintegral power $z^{\rho}$ with $\rho \in \mathbb{R} \setminus
\mathbb{N}_{0}$ is not analytic at the origin. Accordingly, a power
series for $z^{\rho}$ about $z=0$ cannot exist. At least formally, the
\emph{infinite} Laguerre expansion (\ref{GenPow2GLag}) for $z^{\rho}$ can
be rearranged to yield the power series (\ref{Chk_x^m_GlagPol_6}), but
this power series is mathematically meaningless since almost all of its
series coefficients are infinite in magnitude.

Nevertheless, the example of the general power function $z^{\rho}$ is
instructive since it shows that mathematics cannot be cheated by formally
rearranging Laguerre series. Therefore, we can try to use the approach
described above also in the case of essentially arbitrary infinite
Laguerre series of the type of (\ref{f_Exp_GLag}). We insert the explicit
expression (\ref{GLag_1F1}) of the generalized Laguerre polynomial into
the Laguerre series for some function, and rearrange the order of
summations. As a final step, we have to analyze whether and under which
conditions the inner infinite series expansions for the coefficients of
the resulting the power series converge.

But first, let us consider a partial sum of the general Laguerre series
(\ref{f_Exp_GLag}):
\begin{equation}
  \label{FinSum_GLag}
  f_{N} (z) \; = \; \sum_{n=0}^{N} \,
  \lambda_{n}^{(\alpha)} \, L_{n}^{(\alpha)} (z) \, ,
  \qquad N \in \mathbb{N}_{0} \, .
\end{equation}
Such a finite sum is simply a polynomial in $z$ and it is always possible
to reformulate it by interchanging the order of the nested finite
summations. If we insert the explicit expression (\ref{GLag_1F1}) of the
generalized Laguerre polynomials into (\ref{FinSum_GLag}) and rearrange
the order of summations, we obtain after some algebra:
\begin{align}
  f_{N} (z) & \; = \; \sum_{n=0}^{N} \, \lambda_{n}^{(\alpha)} \,
  \frac{(\alpha+1)_{n}}{n!} \, \sum_{\nu=0}^{n} \,
  \frac{(-n)_{\nu}}{(\alpha+1)_{\nu}} \, \frac{z^{\nu}}{\nu!}
  \\
  & \; = \; \sum_{\nu=0}^{N} \, \frac{z^{\nu}}{(\alpha+1)_{\nu} \nu!} \,
  \sum_{n=\nu}^{N} \, \frac{(-n)_{\nu} (\alpha+1)_{n}}{n!} \,
  \lambda_{n}^{(\alpha)} \, .
\end{align}
This expression can be streamlined further, yielding
\begin{equation}
  \label{RearrFinSum_GLag}
  f_{N} (z) \; = \; \sum_{\nu=0}^{N} \, \frac{(-z)^{\nu}}{\nu!} \,
  \sum_{\mu=0}^{N-\nu} \, \frac{(\alpha+\nu+1)_{\mu}}{\mu!} \,
  \lambda_{\mu+\nu}^{(\alpha)} \, .
\end{equation}
Thus, in the case of \emph{finite} Laguerre expansions of the type of
(\ref{FinSum_GLag}), a rearrangement of the order of the nested finite
summations is always possible.

Let us now consider the rearrangement of the infinite series
(\ref{f_Exp_GLag}). By inserting the explicit expression (\ref{GLag_1F1})
of the generalized Laguerre polynomials into (\ref{f_Exp_GLag}) and
rearranging the order of summations, we formally obtain the following
power series in $z$:
\begin{equation}
  \label{Rearr_f_Exp_GLag}
  f (z) \; = \;
  \sum_{\nu=0}^{\infty} \, \frac{(-z)^{\nu}}{\nu!} \,
  \sum_{\mu=0}^{\infty} \, \frac{(\alpha+\nu+1)_{\mu}}{\mu!} \,
  \lambda_{\mu+\nu}^{(\alpha)} \, .
\end{equation}

If we compare (\ref{RearrFinSum_GLag}) and (\ref{Rearr_f_Exp_GLag}), we
immediately see that the rearrangement of an infinite Laguerre expansion
is not necessarily possible since we now have an inner \emph{infinite
  series} instead of an inner \emph{finite sum}. Accordingly, many things
can go wrong if we mechanically perform the limit $N \to \infty$ in
(\ref{RearrFinSum_GLag}). The power series (\ref{Rearr_f_Exp_GLag}) for
$f (z)$ makes sense if and only if the inner series on the right-hand
side of (\ref{Rearr_f_Exp_GLag}) converges for every $\nu \in
\mathbb{N}_{0}$. Otherwise, we have a formal power series with expansion
coefficients that are infinite in magnitude. This scenario corresponds to
the formal, but mathematically meaningless power series
(\ref{Chk_x^m_GlagPol_6}) for $z^{\rho}$ with $\rho \in \mathbb{R}
\setminus \mathbb{N}_{0}$.

\typeout{==> Section: Algebraically Decaying Series Coefficients}
\section{Algebraically Decaying Series Coefficients}
\label{Sec:AlgebraicallyDecayingSeriesCoefficients}

In this Section, the convergence of the inner $\mu$ series in the
rearranged Laguerre expansion (\ref{Rearr_f_Exp_GLag}) is analyzed by
making several assumptions on the large index ($n \to \infty$)
asymptotics of the coefficients $\lambda_{n}^{(\alpha)}$.

The sufficient convergence condition (\ref{ConvCon_LambdaSer}) shows that
the coefficients $\lambda_{n}^{(\alpha)}$ in (\ref{f_Exp_GLag}) can decay
algebraically in $n$, which in practice implies (very) bad convergence.
Thus, let us assume for the moment that the $\lambda_{n}^{(\alpha)}$ all
have the same sign at least for sufficiently large indices $n$, and that
they possess the following large index asymptotics:
\begin{equation}
 \label{Asy_lambda_1}
  \lambda_{n}^{(\alpha)} \; \sim \; n^{-\beta} \, ,
  \qquad n \to \infty \, , \quad \beta > 0 \, .
\end{equation}

For an analysis of the convergence of the inner $\mu$ series, it is
helpful to rewrite (\ref{Rearr_f_Exp_GLag}) as follows:
\begin{equation}
  \label{Rearr_f_Exp_GLag_1}
  f (x) \; = \; \frac{1}{\Gamma (\alpha+1)} \,
  \sum_{\nu=0}^{\infty} \, \frac{(-x)^{\nu}}{(\alpha+1)_{\nu} \, \nu!} \,
  \sum_{\mu=0}^{\infty} \, \frac{\Gamma (\alpha+\mu+\nu+1)}{\mu!} \,
  \lambda_{\mu+\nu}^{(\alpha)} \, .
\end{equation}
We first analyze the large index asymptotics of the factor $\Gamma
(\alpha+\mu+\nu+1)/\mu!$. With the help of (\ref{AsyGammaRatio}), we
obtain for fixed and finite $\nu \in \mathbb{N}_{0}$ the following
asymptotic approximation:
\begin{equation}
  \label{CF_asy_mu_ser}
  \frac{\Gamma (\alpha+\mu+\nu+1)}{\mu!}
  \; \sim \; \mu^{\alpha+\nu} \, , \qquad \mu \to \infty \, .
\end{equation}
For fixed and finite $\nu \in \mathbb{N}_{0}$, the asymptotic estimate
(\ref{Asy_lambda_1}) translates to
\begin{equation}
  \label{Asy_lambda_AlgebraicMonotone}
  \lambda_{\mu+\nu}^{(\alpha)} \; \sim \; (\mu+\nu)^{-\beta} \; = \;
  \mu^{-\beta} \, + \, \mathrm{O} \bigl( \mu^{-\beta-1} \bigr) \, ,
  \qquad \mu \to \infty \, .
\end{equation}
Combination of (\ref{CF_asy_mu_ser}) and
(\ref{Asy_lambda_AlgebraicMonotone}) yields
\begin{equation}
  \label{LargeInd_MuSum_1}
  \frac{\Gamma (\alpha+\mu+\nu+1)}{\mu!} \, \lambda_{\mu+\nu}^{(\alpha)}
  \; \sim \; \mu^{\alpha+\nu-\beta} \, , \qquad \mu \to \infty \, .
\end{equation}
Thus, the inner $\mu$ series in (\ref{Rearr_f_Exp_GLag}) diverges at
least for sufficiently large values of the outer index $\nu$ if the
series coefficients $\lambda_{n}^{(\alpha)}$ occurring in
(\ref{Rearr_f_Exp_GLag}) ultimately have the same sign and decay
algebraically like a fixed power $\beta$ of the index $n$. Thus, a
function represented by a Laguerre expansion with ultimately monotone and
algebraically decaying series coefficients cannot be analytic in a
neighborhood of the origin. This conclusion is in agreement with the
remarks by Gottlieb and Orszag \cite[p.\ 42]{Gottlieb/Orszag/1977} and by
Doha \cite[p.\ 5452]{Doha/2003}, respectively, who had stated that the
Laguerre series for a given function converges faster than algebraically
if the function under consideration is analytic at the origin.

As discussed in more details in Section
\ref{Sec:LagSerGeneralPowerFunction}, the general power function
$z^{\rho}$ with $\rho \in \mathbb{R} \setminus \mathbb{N}_{0}$ is not
analytic at the origin. This fact can also be deduced from its Laguerre
series (\ref{GenPow2GLag}). With the help of (\ref{AsyGammaRatio}), we
obtain the following leading order asymptotic estimate for the
coefficients in (\ref{GenPow2GLag}):
\begin{equation}
  \label{Asy_Cf_x^mu_LagSer}
  \frac{\Gamma (-\rho+n)}{\Gamma (\alpha+n+1)} \; \sim \;
  n^{-\alpha-\rho-1} \, , \qquad n \to \infty \, .
\end{equation}
Comparison with (\ref{ConvCon_LambdaSer}) shows that this asymptotic
estimate implies the convergence of the Laguerre series
(\ref{GenPow2GLag}) for $z^{\rho}$ with respect to the norm
(\ref{Def:Norm_Lag}) of the Hilbert space $L^{2}_{x^{\alpha}
  \mathrm{e}^{-x}} \bigl( [0, \infty) \bigr)$ if $\alpha+2\rho > -1$
holds. However, this estimate also shows that $z^{\rho}$ with $\rho \in
\mathbb{R} \setminus \mathbb{N}_{0}$ cannot be analytic at the origin
$z=0$.

The convergence properties of \emph{monotone} and \emph{alternating}
series differ substantially. A monotone series $\sum_{n=0}^{\infty}
a_{n}$, whose terms all have the same sign, converges, if the series
terms $a_{n}$ decay at least like $a_{n} = \mathrm{O} \bigl(
n^{-1-\epsilon} \bigr)$ with $\epsilon > 0$ as $n \to \infty$ or faster.
In contrast, an alternating series $\sum_{n=0}^{\infty} (-1)^{n} \vert
b_{n} \vert$ converges if the terms $b_{n}$ decrease in magnitude and
approach zero as $n \to \infty$. This alone would not suffice to
guarantee the convergence of the inner $\mu$ series in
(\ref{Rearr_f_Exp_GLag}). However, alternating series have the undeniable
advantage that summability techniques for divergent series can be
employed (further details as well as numerous references can be found in
Appendices \ref{App:DivergentSeries} and
\ref{App:SequenceTransformations}). In this way, it is frequently
possible to associate a finite value to a divergent alternating series
$\sum_{n=0}^{\infty} (-1)^{n} \vert b_{n} \vert$, whose terms do not
vanish as $n \to \infty$ and even grow in magnitude with increasing
index.

Thus, it should be interesting to investigate whether Laguerre series of
the type of (\ref{f_Exp_GLag}) with \emph{strictly alternating} series
coefficients $\lambda_{n}^{(\alpha)}$, whose absolute values are at least
asymptotically proportional to a fixed power of the index $n$, correspond
to functions that are analytic in a neighborhood of the origin.

As an example of a Laguerre series with ultimately strictly alternating
and algebraically decaying coefficients, let us consider the following
Laguerre series, which differs from the Laguerre series
(\ref{GenPow2GLag}) for the general power function $z^{\rho}$ only by the
alternating sign $(-1)^{n}$:
\begin{align}
  \label{LagSer_AltSerCoeffs}
  G_{\rho}^{(\alpha)} (z) & \; = \;
  \frac{\Gamma (\rho+\alpha+1)}{\Gamma (\alpha+1} \,
  \sum_{n=0}^{\infty} \, (-1)^{n} \,
  \frac{(-\rho)_{n}}{(\alpha+1)_{n}}
  \, L_{n}^{(\alpha)} (z) \, ,
  \notag \\
  & \qquad \rho \in \mathbb{R} \setminus \mathbb{N}_0 \, ,
  \qquad \alpha+2\rho > - 1 \, .
\end{align}
The convergence condition $\alpha+2\rho > - 1$ guarantees that
$G_{\rho}^{(\alpha)} (z)$ belongs just like the power function $z^{\rho}$
with $\rho \in \mathbb{R} \setminus \mathbb{N}_{0}$ to the Hilbert space
$L^{2}_{x^{\alpha} \mathrm{e}^{-x}} \bigl( [0, \infty) \bigr)$. The case
$\rho \in \mathbb{N}_{0}$ is excluded because then $G_{\rho}^{(\alpha)}$
would be a polynomial in $z$, whose analyticity at the origin is obvious.

It is, however, unclear whether $G_{\rho}^{(\alpha)} (z)$ with $\rho \in
\mathbb{R} \setminus \mathbb{N}_{0}$ is analytic at the origin, i.e.,
whether the inner $\mu$ series in the rearranged Laguerre series
(\ref{Rearr_f_Exp_GLag}) converges for arbitrary values of $\nu \in
\mathbb{N}_{0}$. For an investigation of this question, let us define
\begin{equation}
  \label{Def_LanbdaG}
  \lambda_{n}^{(\alpha)} \; = \; (-1)^{n} \,
  \frac{\Gamma (\rho+\alpha+1)}{\Gamma (-\rho)} \,
  \frac{\Gamma (-\rho+n)}{\Gamma (\alpha+n+1)} \, .
\end{equation}
Inserting this into the modified rearranged Laguerre expansion
(\ref{Rearr_f_Exp_GLag_1}) yields:
\begin{align}
  \label{Rearr_G_Exp_GLag_1}
  G_{\rho}^{(\alpha)} (z) & \; = \;
  \frac{\Gamma (\rho+\alpha+1)}{\Gamma (\alpha+1)} \,
  \sum_{\nu=0}^{\infty} \, \frac{(-\rho)_{\nu}}{(\alpha+1)_{\nu}} \,
  \frac{z^{\nu}}{\nu!} \, \sum_{\mu=0}^{\infty} \, (-1)^{\mu} \,
  \frac{(-\rho+\nu)_{\mu}}{\mu!}
  \\
  \label{Rearr_G_Exp_GLag_3}
  & \; = \;
  \frac{\Gamma (\rho+\alpha+1)}{\Gamma (\alpha+1)} \,
  \sum_{\nu=0}^{\infty} \, \frac{(-\rho)_{\nu}}{(\alpha+1)_{\nu}} \,
  \frac{z^{\nu}}{\nu!} \,
  {}_{1} F_{0} \bigl( -\rho+\nu; -1 \bigr) \, .
\end{align}
The generalized hypergeometric series ${}_{1} F_{0}$ in
(\ref{Rearr_G_Exp_GLag_3}) is a special case of the binomial series
(\ref{BinomSer}) which can be expressed in close form. Thus, we obtain:
\begin{equation}
  \label{1F0_Res_1}
  {}_{1} F_{0} \bigl( -\rho+\nu; -1 \bigr) \; = \;
  \lim_{z \to -1} \, (1-z)^{\rho-\nu} \; = \; 2^{\rho-\nu} \, .
\end{equation}
Since (\ref{AsyGammaRatio}) implies $(-\rho+\nu)_{\mu}/\mu! = \mathrm{O}
\bigl( \mu^{\nu-\rho-1} \bigr)$ as $\mu \to \infty$, the hypergeometric
series ${}_{1} F_{0} \bigl( -\rho+\nu; z \bigr)$ converges only for
$\vert z \vert < 1$. Thus, ${}_{1} F_{0} \bigl( -\rho+\nu; -1 \bigr)$ is
-- strictly speaking -- undefined and a divergent series. However,
$(1-z)^{n-\rho}$ remains well defined as $z \to -1$.  Therefore,
(\ref{1F0_Res_1}) essentially corresponds to an analytic continuation
that implicitly uses the concepts of \emph{Abel summation} which is for
instance discussed in Hardy's classic book \cite{Hardy/1949}.

Inserting (\ref{1F0_Res_1}) into (\ref{Rearr_G_Exp_GLag_3}) yields:
\begin{equation}
  \label{Rearr_G_Exp_GLag_4}
  G_{\rho}^{(\alpha)} (z) \; = \; 2^{\rho} \,
  \frac{\Gamma (\rho+\alpha+1)}{\Gamma (\alpha+1)} \,
  {}_{1} F_{1} \bigl( -\rho; \alpha+1; z/2 \bigr) \, .
\end{equation}
A confluent hypergeometric series ${}_{1} F_{1} (a; b; z)$ converges
absolutely for all complex $a$, $b$, and $z$ as long as $-b \notin
\mathbb{N}_{0}$ (see for example \cite[p.\ 2]{Slater/1960}). Thus, for
$-b \notin \mathbb{N}_{0}$ such a ${}_{1} F_{1}$ is an \emph{entire}
function in $z$. Since we always assume $\alpha > -1$,
$G_{\rho}^{(\alpha)} (z)$ is in every neighborhood of the origin $z=0$ an
\emph{analytic} function. Thus, (\ref{Rearr_G_Exp_GLag_4}) shows that the
remarks by Gottlieb and Orszag \cite[p.\ 42]{Gottlieb/Orszag/1977} and by
Doha \cite[p.\ 5452]{Doha/2003}, who stated that such a Laguerre series
converges faster than algebraically if the function under consideration
is analytic at the origin, is imprecise.

If the argument $z$ of a confluent hypergeometric function ${}_{1} F_{1}
(a; b, z)$ is real and approaches $+\infty$, then we have the following
asymptotic behavior (see for example \cite[Eq.\ (4.1.7)]{Slater/1960}):
\begin{equation}
  \label{Asy_1F1_x->=Inf}
  {}_{1} F_{1} (a; b; z) \; = \; \frac{\Gamma (a)}{\Gamma (b)} \,
  \mathrm{e}^{z} \, z^{a-b} \, \bigl[ 1 + \mathrm{O} (1/z) \bigr] \, ,
  \qquad z \to +\infty \, .
\end{equation}
This asymptotic estimate shows that the integral
\begin{align}
  \label{}
  & \int_{0}^{\infty} \, \mathrm{e}^{-z} \, z^{\alpha}
  \bigl[ G_{\rho}^{(\alpha)} (z) \bigr]^{2}
  \mathrm{d} z
  \notag \\
  & \qquad \; = \; \left[ 2^{\rho} \,
  \frac{\Gamma (\rho+\alpha+1)}{\Gamma (\alpha+1)} \right]^{2} \,
  \int_{0}^{\infty} \, \mathrm{e}^{-z} \, z^{\alpha}
  \bigl[ {}_{1} F_{0} \bigl( -\rho; \alpha+1; z/2 \bigr) \bigr]^{2}
  \mathrm{d} z
\end{align}
converges if $\alpha > -1$ and $\alpha+2\rho > -1$ hold. At the lower
integration limit $z=0$, the integrand behaves like $z^{\alpha}$, which
requires $\alpha > -1$, and at the upper integration limit $z=\infty$,
the integrand behaves like $z^{-\alpha-2\rho-2}$, which requires
$\alpha+2\rho > -1$.  Thus, $G_{\rho}^{(\alpha)} (z)$ belongs for $\alpha
> -1$ and $\alpha+2\rho > -1$ to the Hilbert space $L^{2}_{z^{\alpha}
  \mathrm{e}^{-z}} \bigl( [0, \infty) \bigr)$ defined by
(\ref{HilbertL^2_Lag}).

It is possible to check the correctness of the power series
representation (\ref{Rearr_G_Exp_GLag_4}) for $G_{\rho}^{(\alpha)} (z)$
by expanding it in terms of generalized Laguerre polynomials. According
to (\ref{f_Exp_GLag_b}), we then have to compute the following inner
product:
\begin{align}
  \label{lambda_G_2}
  \lambda_{n}^{(\alpha)} & \; = \;
  2^{\rho} \, \frac{n!}{\Gamma (\alpha+n+1)} \,
  \frac{\Gamma (\rho+\alpha+1)}{\Gamma (\alpha+1)}
  \notag \\
  & \qquad \times \, \int_{0}^{\infty} \, \mathrm{e}^{-z} \, z^{\alpha} \,
  L_{n}^{(\alpha)} (z) \,
  {}_{1} F_{1} \bigl( -\rho; \alpha+1; z/2 \bigr) \, \mathrm{d} z \, .
\end{align}
Since the confluent hypergeometric function ${}_{1} F_{1}$ in
(\ref{lambda_G_2}) is analytic in every neighborhood of the origin,
integration and summation can be interchanged. We now use \cite[Eq.\
(7.414.11) on p.\ 850]{Gradshteyn/Rhyzhik/1994}
\begin{equation}
  \label{GR_7.414.11}
  \int_{0}^{\infty} \, \mathrm{e}^{-t} \, t^{\gamma-1} \,
  L_{n}^{(\mu)} (t) \, \mathrm{d} t \; = \;
  \frac{\Gamma (\gamma) \Gamma (1+\mu-\gamma+n)}{n! \Gamma (1+\mu-\gamma)}
  \, , \qquad \Re (\gamma) > 0 \, ,
\end{equation}
to obtain:
\begin{align}
  \label{G_GLag_Int_2}
  & \int_{0}^{\infty} \, \emph{e}^{-z} \, z^{\alpha} \, L_{n}^{(\alpha)}
  (z) \, {}_{1} F_{1} \bigl( -\rho; \alpha+1; z/2 \bigr) \, \mathrm{d} z
  \notag \\
  & \qquad \; = \; \frac{\Gamma (\alpha+1)}{n!} \, \sum_{m=0}^{\infty} \,
  \frac{(-\rho)_{m} \, (-m)_{n}} {m!} \, 2^{-m} \, .
\end{align}
The fact that the Pochhammer symbol $(-m)_{n}$ with $m, n \in
\mathbb{N}_{0}$ satisfies $(-m)_{n} = 0$ for $m = 0, 1, \dots, n-1$
suggest the substitution $m \to n + \nu$, yielding:
\begin{align}
  \label{G_GLag_Int_5}
  & \int_{0}^{\infty} \, \emph{e}^{-z} \, z^{\alpha} \, L_{n}^{(\alpha)}
  (z) \, {}_{1} F_{1} \bigl( -\rho; \alpha+1; z/2 \bigr) \, \mathrm{d} z
  \notag \\
  & \qquad \; = \; (-1)^{n} \, \frac{\Gamma (\alpha+1) \,
    (-\rho)_{n}}{2^{n} \, n!} \, {}_{1} F_{0} \bigl( -\rho+n; 1/2 \bigr)
  \, .
\end{align}
With the help of (\ref{BinomSer}) we obtain
\begin{equation}
  \label{1F1_G_GLag_Int_5}
  {}_{1} F_{0} \bigl( -\rho+n; 1/2 \bigr) \; = \; (1/2)^{\rho-n}
  \; = \; 2^{n-\rho} \, .
\end{equation}
Combination of (\ref{G_GLag_Int_5}) and (\ref{1F1_G_GLag_Int_5}) then
yields:
\begin{align}
  \label{G_GLag_Int_6}
  & \int_{0}^{\infty} \, \emph{e}^{-z} \, z^{\alpha} \,
  L_{n}^{(\alpha)} (z) \,
  {}_{1} F_{1} \bigl( -\rho; \alpha+1; z/2 \bigr) \, \mathrm{d} z
  \notag  \\
  & \qquad \; = \; \frac{(-1)^{n}}{2^{\rho}} \,
  \frac{\Gamma (\alpha+1) \, (-\rho)_{n}}{n!} \, .
\end{align}
We also obtain this result if we combine (\ref{Def_LanbdaG}) and
(\ref{lambda_G_2}).

As a more complicated example, let us now consider the following Laguerre
series:
\begin{equation}
  \label{LagSer_2H1_a_b_c_alpha}
  \prescript{}{2}{\mathcal{H}}_{1}^{(\alpha)} (a, b; c; z) \; = \;
  \sum_{n=0}^{\infty} \, (-1)^{n} \,
  \frac{(a)_{n} (b)_{n}}{(c)_{n} (\alpha+1)_{n}}
  \, L_{n}^{(\alpha)} (z) \, .
\end{equation}
Concerning the real parameters $a$, $b$, and $c$, we assume for the
moment only that this series exists, which requires $-c \notin
\mathbb{N}_{0}$, that it does not terminate, which requires $-a, -b
\notin \mathbb{N}_{0}$. and that it converges with respect to the norm
(\ref{Def:Norm_Lag}) of the Hilbert space $L^{2}_{z^{\alpha}
  \mathrm{e}^{-z}} \bigl( [0, \infty) \bigr)$. It follows from
(\ref{AsyGammaRatio}) that this is the case if $c-a-b+(\alpha+1)/2 > 0$.

For an investigation of the analyticity of
$\prescript{}{2}{\mathcal{H}}_{1}^{(\alpha)} (a, b; c; z)$ at the origin,
let us define
\begin{equation}
  \label{Def_LanbdaH}
  \lambda_{n}^{(\alpha)} \; = \; (-1)^{n} \,
  \frac{(a)_{n} (b)_{n}}{(c)_{n} (\alpha+1)_{n}} \, .
\end{equation}
By proceeding as in the case of the Laguerre series
(\ref{LagSer_AltSerCoeffs}) for $G_{\rho}^{(\alpha)} (z)$ we obtain after
a short calculation:
\begin{equation}
  \label{Rearr_2H1_Exp_GLag_3}
  \prescript{}{2}{\mathcal{H}}_{1}^{(\alpha)} (a, b; c; z) \; = \;
  \sum_{\nu=0}^{\infty} \, {}_{2} F_{1} ( a+\nu, b+\nu;
  c+\nu; -1) \, \frac{(a)_{\nu} \, (b)_{\nu}}{(c)_{\nu} \, \nu!} \, \,
  \frac{z^{\nu}}{(\alpha+1)_{\nu}} \, .
\end{equation}
It follows from (\ref{AsyGammaRatio}) that $\Gamma (a+\nu+n) \Gamma
(b+\nu+n)/[\Gamma (c+\nu+n) n!] \sim n^{a+b+\nu-c-1}$ as $n \to \infty$.
Consequently, the Gaussian hypergeometric series ${}_{2} F_{1} \bigl(
a+\nu, b+\nu; c+\nu; z \bigr)$ with $z=-1$ diverges for all sufficiently
large $\nu \in \mathbb{N}_{0}$. Nevertheless, it is possible to associate
a finite value to the hypergeometric function corresponding to this
divergent series with the help of the following two linear
transformations \cite[p.\ 47]{Magnus/Oberhettinger/Soni/1966}:
\begin{align}
  \label{LTr_1}
  {}_2 F_1 (a, b; c; z) & \; = \; (1-z)^{-a} \, {}_2 F_1 \bigl( a, c-b;
  c; z/(z-1) \bigr)
  \\
  \label{LTr_2}
  & \; = \; (1-z)^{-b} \, {}_2 F_1 \bigl( c-a, b; c; z/(z-1) \bigr) \, .
\end{align}
The transformation $z \to z' = z/(z-1)$ maps $z=-1$, which is located on
the boundary of the circle of convergence, to $z' = 1/2$, which is
located in the interior of the circle of convergence. Thus, the in
general \emph{divergent} hypergeometric series in
(\ref{Rearr_2H1_Exp_GLag_3}) can be replaced by a \emph{convergent}
hypergeometric series according to
\begin{align}
  \label{Rearr_2H1_Exp_GLag_3_LTr_1}
  {}_2 F_1 (a+\nu, b+\nu; c+\nu; -1) & \; = \; 2^{-a-\nu} \, {}_2 F_1
  \bigl( a+\nu, c-b; c+\nu; 1/2 \bigr)
  \\
  \label{Rearr_2H1_Exp_GLag_3_LTr_2}
  & \; = \; 2^{-b-\nu} \, {}_2 F_1 \bigl( c-a, b+\nu; c+\nu; 1/2 \bigr)
  \, .
\end{align}
Inserting (\ref{Rearr_2H1_Exp_GLag_3_LTr_1}) and
(\ref{Rearr_2H1_Exp_GLag_3_LTr_2}) into (\ref{Rearr_2H1_Exp_GLag_3})
yields two equivalent power series expansions for the function defined by
the Laguerre series (\ref{LagSer_2H1_a_b_c_alpha}), which seem to be new:
\begin{align}
  \label{Rearr_2H1_Exp_GLag_4}
  \prescript{}{2}{\mathcal{H}}_{1}^{(\alpha)} (a, b; c; z)
  & \; = \; 2^{-a} \, \sum_{\nu=0}^{\infty} \, {}_2 F_1 \bigl(
  a+\nu, c-b; c+\nu; 1/2 \bigr) \, \frac{(a)_{\nu} \,
    (b)_{\nu}}{(c)_{\nu} \, \nu!} \, \,
  \frac{(z/2)^{\nu}}{(\alpha+1)_{\nu}}
  \\
  \label{Rearr_2H1_Exp_GLag_5}
  & \; = \; 2^{-b} \, \sum_{\nu=0}^{\infty} \, {}_2 F_1 \bigl(
  c-a, b+\nu; c+\nu; 1/2 \bigr) \, \frac{(a)_{\nu} \,
    (b)_{\nu}}{(c)_{\nu} \, \nu!} \, \,
  \frac{(z/2)^{\nu}}{(\alpha+1)_{\nu}} \, .
\end{align}
A detailed analysis of the domain of analyticity of
$\prescript{}{2}{\mathcal{H}}_{1}^{(\alpha)} (a, b; c; z)$ requires
asymptotic estimates of the behavior of the Gaussian hypergeometric
series in (\ref{Rearr_2H1_Exp_GLag_3}), (\ref{Rearr_2H1_Exp_GLag_4}) and
(\ref{Rearr_2H1_Exp_GLag_5}) as $\nu \to \infty$. If we use the linear
transformation \cite[p.\ 47]{Magnus/Oberhettinger/Soni/1966}
\begin{equation}
  \label{LTr_0}
  {}_2 F_1 (a, b; c; z) \; = \;
  (1-z)^{c-a-b} \, {}_2 F_1 (c-a, c-b; c; z) \, ,
\end{equation}
which certainly holds for $\vert z \vert < 1$, we obtain:
\begin{equation}
  \label{AnCon_Rearr_2H1_Exp_GLag_2F1_1}
  {}_2 F_1 (a+\nu, b+\nu; c+\nu; z) \; = \;
  (1-z)^{c-a-b-\nu} \, {}_2 F_1 (c-a, c-b; c+\nu; z) \, .
\end{equation}
The asymptotic estimate (\ref{AsyGammaRatio}) yields $\Gamma (c-a+n)
\Gamma (c-b+n)/[\Gamma (c+\nu+n) n!] \sim n^{c-a-b-\nu-1}$ as $n \to
\infty$. Accordingly, the hypergeometric series on the right-hand side of
(\ref{AnCon_Rearr_2H1_Exp_GLag_2F1_1}) converges at least for
sufficiently large values of $\nu$ also for $z=-1$ and provides an
analytic continuation. Moreover, (\ref{AnCon_Rearr_2H1_Exp_GLag_2F1_1})
is a convenient starting point for the construction of an asymptotic
expansion of the divergent hypergeometric series in
(\ref{Rearr_2H1_Exp_GLag_3}) as $\nu \to \infty$ (compare also \cite[Eq.\
(15)]{Temme/2003}):
\begin{align}
  \label{Asy_ny_large_2F1_1}
  & {}_2 F_1 (a+\nu, b+\nu; c+\nu; -1) \; = \;
  2^{c-a-b-\nu} \, \sum_{n=0}^{\infty} \, (-1)^{n} \,
  \frac{(c-a)_{n} (c-b)_{n}}{(c+\nu)_{n} n!}
  \\
  \label{Asy_ny_large_2F1_2}
  & \qquad \; = \; 2^{c-a-b-\nu} \,
  \left[ 1 - \frac{(c-a) (c-b)}{(c+\nu)} +
  \mathrm{O} \bigl( \nu^{-2} \bigr) \right] \, ,
  \qquad \nu \to \infty \, .
\end{align}
This asymptotic estimate shows that
$\prescript{}{2}{\mathcal{H}}_{1}^{(\alpha)} (a, b; c; z)$ is for $-c
\notin \mathbb{N}_{0}$ analytic in every neighborhood of the origin.

\typeout{==> Section: Exponentially and Factorially Decaying Series
  Coefficients}
\section{Exponentially and Factorially Decaying Series Coefficients}
\label{Sec:ExponentiallyFactoriallyDecayingSeriesCoefficients}

It is immediately obvious that the inner $\mu$ series in
(\ref{Rearr_f_Exp_GLag}) converges if the series coefficients
$\lambda_{n}^{(\alpha)}$ decay for sufficiently large indices $n$
exponentially, satisfying for instance
\begin{equation}
  \label{ExpoDecay_lambda}
  \lambda_{n}^{(\alpha)} \; \sim \; n^{\theta} \, R^{n} \, ,
  \qquad \theta \in \mathbb{R} \, , \qquad \vert R \vert < 1 \, ,
  \qquad n \to \infty \, .
\end{equation}

As the probably most simple example of a Laguerre series of the type of
(\ref{f_Exp_GLag}) with exponentially decaying coefficients, let us
consider the following expansion:
\begin{equation}
  \label{calE_LagSer}
  \mathcal{E}^{(\alpha)} (t; z) \; = \;
  \sum_{n=0}^{\infty} \, t^{n} \, L_{n}^{(\alpha)} (z) \, .
\end{equation}
It is immediately obvious that the Laguerre series (\ref{calE_LagSer})
converges in the mean with respect to the norm (\ref{Def:Norm_Lag}) of
the Hilbert space $L^{2}_{z^{\alpha} \mathrm{e}^{-z}} \bigl( [0, \infty)
\bigr)$ for $\vert t \vert < 1$, and that it diverges for $\vert t \vert
\ge 1$. Equivalently, we can say that the series (\ref{MeanConvCond_f})
for the norm $\bigl\Vert \mathcal{E}^{(\alpha)} (t; z)
\bigr\Vert_{z^{\alpha} \mathrm{e}^{-z}, 2}$ produces a finite result for
$\vert t \vert < 1$ and that it diverges for $\vert t \vert \ge 1$.

Ignoring for the moment all questions of convergence, let us now assume
that $t$ is an unspecified complex number and insert
$\lambda_{n}^{(\alpha)} = t^{n}$ into (\ref{Rearr_f_Exp_GLag}). Then, we
obtain:
\begin{align}
  \label{calE_PowSer_1}
  \mathcal{E}^{(\alpha)} (t; z) & \; = \;
  \sum_{\nu=0}^{\infty} \, \frac{(-t z)^{\nu}}{\nu!} \,
  \sum_{\mu=0}^{\infty} \, \frac{(\alpha+\nu+1)_{\mu}}{\mu!} \, t^{\mu}
  \\
  \label{calE_PowSer_2}
  & \; = \;
  \sum_{\nu=0}^{\infty} \, \frac{(-t z)^{\nu}}{\nu!} \,
  {}_{1} F_{0} \bigl( \alpha+\nu+1; t \bigr) \, .
\end{align}
A nonterminating hypergeometric series ${}_{1} F_{0} (a; z)$ converges
only in the interior of the unit circle. Since we always assume $\alpha >
-1$, the hypergeometric series ${}_{1} F_{0} \bigl( \alpha+\nu+1; t
\bigr)$ in (\ref{calE_PowSer_2}) does not terminate and we have to
require $\vert t \vert < 1$. However, ${}_{1} F_{0} \bigl( \alpha+\nu+1;
t \bigr)$ can according to (\ref{BinomSer}) be expressed in closed form,
\begin{equation}
  \label{1F0_calE_PowSer}
  {}_{1} F_{0} \bigl( \alpha+\nu+1; t \bigr) \; = \;
  (1-t)^{-\alpha-\nu-1} \, ,
  \qquad t \in \mathbb{C} \setminus \{ 1 \} \, ,
\end{equation}
which yields an analytic continuation for $t \ne 1$.  Inserting
(\ref{1F0_calE_PowSer}) into (\ref{calE_PowSer_2}) yields:
\begin{align}
  \label{calE_PowSer_3}
  \mathcal{E}^{(\alpha)} (t; z) & \; = \; (1-t)^{-\alpha-1} \,
  \sum_{\nu=0}^{\infty} \, \frac{\bigl(-t z/[1-t] \bigr)^{\nu}}{\nu!}
  \\
  \label{calE_PowSer_4}
  & \; = \; (1-t)^{-\alpha-1} \, \exp \bigl(-t z/[1-t] \bigr) \, .
\end{align}
Accordingly, the Laguerre series (\ref{calE_LagSer}) for
$\mathcal{E}^{(\alpha)\ } (t; z)$ is nothing but the well known
generating function \cite[p.\ 242]{Magnus/Oberhettinger/Soni/1966}
\begin{equation}
  \label{GLag_Gen_1}
  \sum_{n=0}^{\infty} \, L_{n}^{(\alpha)} (z) \, t^n \; = \;
  (1-t)^{-\alpha-1} \, \exp \bigl( zt/[t-1] \bigr) \, , \qquad \vert t
  \vert < 1 \, ,
\end{equation}
in disguise.

As a generalization of the alternating Laguerre series
(\ref{LagSer_AltSerCoeffs}), let us consider the following expansion with
exponentially decaying coefficients:
\begin{align}
  \label{LagSer_ExpoDecaySerCoeffs}
  \mathcal{G}_{\rho}^{(\alpha)} (s; z) & \; = \;
  \frac{\Gamma (\rho+\alpha+1)}{\Gamma (\alpha+1)} \,
  \sum_{n=0}^{\infty} \, (-s)^{n} \,
  \frac{(-\rho)_{n}}{(\alpha+1)_{n}}
  \, L_{n}^{(\alpha)} (z) \, ,
  \notag \\
  & \qquad \rho \in \mathbb{R} \setminus \mathbb{N}_0 \, ,
  \qquad \alpha+2\rho > - 1 \, .
\end{align}
Obviously, this Laguerre series converges in the mean for $\vert s \vert
< 1$.

For an investigation of the analyticity of $\mathcal{G}_{\rho}^{(\alpha)}
(s; z)$ at the origin, let us define
\begin{equation}
  \label{Def_LanbdaG_Expo}
  \lambda_{n}^{(\alpha)} \; = \;
  \frac{\Gamma (\rho+\alpha+1)}{\Gamma (-\rho)} \,
  \frac{\Gamma (-\rho+n)}{\Gamma (\alpha+n+1)} \, (-s)^{n} \, .
\end{equation}
Inserting this into the modified rearranged Laguerre expansion
(\ref{Rearr_f_Exp_GLag_1}) yields:
\begin{align}
  \label{Rearr_calG_Exp_GLag_2}
  \mathcal{G}_{\rho}^{(\alpha)} (s; z) & \; = \;
  \frac{\Gamma (\rho+\alpha+1)}{\Gamma (\alpha+1)} \,
  \sum_{\nu=0}^{\infty} \, \frac{(-\rho)_{\nu}}{(\alpha+1)_{\nu}} \,
  \frac{(s z)^{\nu}}{\nu!} \, \sum_{\mu=0}^{\infty} \,
  \frac{(-\rho+\nu)_{\mu}}{\mu!} \, (-s)^{\mu}
  \\
  \label{Rearr_calG_Exp_GLag_3}
  & \; = \;
  \frac{\Gamma (\rho+\alpha+1)}{\Gamma (\alpha+1)} \,
  \sum_{\nu=0}^{\infty} \, \frac{(-\rho)_{\nu}}{(\alpha+1)_{\nu}} \,
  \frac{(s z)^{\nu}}{\nu!} \,
  {}_{1} F_{0} \bigl( -\rho+\nu; -s \bigr) \, .
\end{align}
With the help of (\ref{BinomSer}), we obtain
\begin{equation}
  \label{1F0_Res_2}
  {}_{1} F_{0} \bigl( -\rho+\nu; -s \bigr) \; = \;
  (1+s)^{\rho-\nu} \, , \qquad s \ne -1 \, .
\end{equation}
Inserting (\ref{1F0_Res_2}) into (\ref{Rearr_calG_Exp_GLag_3}) yields:
\begin{equation}
  \label{Rearr_calG_Exp_GLag_4}
  \mathcal{G}_{\rho}^{(\alpha)} (s; z) \; = \; (1+s)^{\rho} \,
  \frac{\Gamma (\rho+\alpha+1)}{\Gamma (\alpha+1)} \,
  {}_{1} F_{1} \bigl( -\rho; \alpha+1; s z/[1+s] \bigr) \, .
\end{equation}
The argument of the confluent hypergeometric series ${}_{1} F_{1}$
becomes infinite for $s=-1$. If we assume $\alpha > -1$, then the
right-hand side of (\ref{Rearr_calG_Exp_GLag_4}) is mathematically
meaningful for all $z \in \mathbb{C}$ and for all $s \in \mathbb{C}
\setminus \{ -1 \}$. However, this is not true for the corresponding
Laguerre series (\ref{LagSer_ExpoDecaySerCoeffs}). The convergence of the
series (\ref{MeanConvCond_f}) for the norm $\bigl\Vert
\mathcal{G}_{\rho}^{(\alpha)} (s; z) \bigr\Vert_{z^{\alpha}
  \mathrm{e}^{-z}, 2}$ is only guaranteed for $\vert s \vert < 1$.

It follows from (\ref{Rearr_calG_Exp_GLag_4}) that
$\mathcal{G}_{\rho}^{(\alpha)} (s; z)$ is essentially a confluent
hypergeometric series ${}_{1} F_{1}$. Accordingly, the Laguerre series
(\ref{LagSer_ExpoDecaySerCoeffs}) for $\mathcal{G}_{\rho}^{(\alpha)} (s;
z)$ can be rewritten as follows:
\begin{equation}
  \label{1F1_LagSer}
  {}_{1} F_{1} \bigl( -\rho; \alpha+1; s z/[1+s] \bigr) \; = \;
  (1+s)^{-\rho} \, \sum_{n=0}^{\infty} \, (-s)^{n} \,
  \frac{(-\rho)_{n}}{(\alpha+1)_{n}} \, L_{n}^{(\alpha)} (z) \, .
\end{equation}
If we make in (\ref{1F1_LagSer}) the substitutions $\rho \to -c$ and $s
\to -t$, we see that the function $\mathcal{G}_{\rho}^{(\alpha)} (s; z)$
is nothing but the well known generating function \cite[Eq.\ (3) on p.\
202]{Rainville/1971})
\begin{equation}
  \label{GLag_Gen_4}
  \sum_{n=0}^{\infty} \, \frac{ (c)_n L_{n}^{(\alpha)}
    (x)}{(\alpha+1)_n} \, t^n \; = \; (1-t)^{-c} \,
  {}_1 F_1 \bigl( c; \alpha+1; -xt/[1-t] \bigr)
\end{equation}
in disguise.

If we set $s=-1$ in the Laguerre series (\ref{LagSer_ExpoDecaySerCoeffs})
for $\mathcal{G}_{\rho}^{(\alpha)} (s; z)$, we obtain the Laguerre series
(\ref{GenPow2GLag}) for $z^{\rho}$. To analyze the behavior of
$\mathcal{G}_{\rho}^{(\alpha)} (s; z)$ as $s \downarrow -1$, we write $s
= -1 + \delta$ with $\delta \ge 0$ and consider in
(\ref{Rearr_calG_Exp_GLag_4}) the limit $\delta \downarrow 0$. We obtain
\begin{equation}
  \label{Rearr_calG_Exp_GLag_5}
  \mathcal{G}_{\rho}^{(\alpha)} (\delta-1; z) \; = \;
  \delta^{\rho} \,
  \frac{\Gamma (\rho+\alpha+1)}{\Gamma (\alpha+1)} \, {}_{1} F_{1}
  \bigl( -\rho; \alpha+1; (\delta-1) z/\delta \bigr) \, .
\end{equation}
Next, we use the following asymptotic estimate \cite[Eq.\
(4.1.8)]{Slater/1960}:
\begin{equation}
  \label{Asy_1F1_z_NegInfArg}
  {}_{1} F_{1} (a; b; z) \; = \; \frac{\Gamma (b)}{\Gamma (b-a)} \,
  (-z)^{-a} \, \bigl[ 1 + \mathrm{O} (1/z) \bigr] \, ,
  \qquad z \to - \infty \, .
\end{equation}
If we insert (\ref{Asy_1F1_z_NegInfArg}) into
(\ref{Rearr_calG_Exp_GLag_5}), we obtain in the limit of vanishing
$\delta \ge 0$:
\begin{align}
  \label{Rearr_calG_Exp_GLag_6_lim}
  \mathcal{G}_{\rho}^{(\alpha)} (-1; z) & \; = \;
  \lim_{\delta \downarrow 0} \,
  \mathcal{G}_{\rho}^{(\alpha)} (\delta-1; z)
  \notag \\
  & \; = \; \lim_{\delta \downarrow 0} \, [(1-\delta)z]^{\rho} \,
  \biggl\{ 1 + \mathrm{O}
  \left( \frac{\delta}{(\delta-1)z} \right) \biggr\}
  \; = \; z^{\rho} \, .
\end{align}
Thus, $\mathcal{G}_{\rho}^{(\alpha)} (s; z)$ possesses the
\emph{one-sided} limit $z^{\rho}$ as $s \downarrow -1$, as it should
according to the Laguerre series (\ref{GenPow2GLag}) and
(\ref{LagSer_ExpoDecaySerCoeffs}).

The transformation formula (\ref{Rearr_f_Exp_GLag}) can only be used if
the function represented by the Laguerre series is analytic at the
origin, because otherwise the inner $\mu$ series diverges. Thus, it is
not possible to obtain in this way an explicit expression for a
nonanalytic function defined by a Laguerre series of the type of
(\ref{f_Exp_GLag}) with monotone and algebraically decaying series
coefficients $\lambda_{n}^{(\alpha)}$.

However, the fact that the in general nonanalytic power function
$z^{\rho}$ can according to (\ref{Rearr_calG_Exp_GLag_6_lim}) be obtained
by considering the one-sided limit $s \downarrow -1$ in the analytic
function $\mathcal{G}_{\rho}^{(\alpha)} (s; z)$ indicates that the
situation is not as hopeless as it may look at first sight.

Let us therefore assume for the moment that the coefficients
$\lambda_{n}^{(\alpha)}$ of a Laguerre series are monotone and decay
algebraically as $n \to \infty$. This implies that the inner $\mu$ series
in (\ref{Rearr_f_Exp_GLag}) produce infinities. However, we could try to
apply the transformation formula (\ref{Rearr_f_Exp_GLag}) to a modified
Laguerre series with coefficients $t ^{n} \lambda_{n}^{(\alpha)}$. Since
these coefficients decay exponentially for $\vert t \vert < 1$, this
modified Laguerre series represents an analytic function. If we succeed
in constructing an explicit expression for this analytic function, then
we can try to perform the one-sided limit $t \uparrow 1$ in this
expression. If this can be done, we obtain an explicit expression for a
nonanalytic function.  Obviously, this idea deserves to be investigated
further.

Although the explicit expression (\ref{Rearr_calG_Exp_GLag_4}) for
$\mathcal{G}_{\rho}^{(\alpha)} (s; z)$ is mathematically well defined for
all $s \in \mathbb{C} \setminus \{ -1 \}$ and possesses a one-sided limit
for $s \downarrow -1$, its Laguerre series
(\ref{LagSer_ExpoDecaySerCoeffs}) requires more restrictive conditions.
In general, generating functions of the kind of
(\ref{LagSer_ExpoDecaySerCoeffs}) converge only for $\vert s \vert < 1$,
i.e., in the interior of the unit circle, and possibly also for
\emph{some} points on the boundary of the unit circle.

If we make in (\ref{LagSer_ExpoDecaySerCoeffs}) or in (\ref{1F1_LagSer})
the substitution $s \to 1$, we obtain the series expansions
(\ref{LagSer_AltSerCoeffs}) and (\ref{Rearr_G_Exp_GLag_4}), respectively,
for $G_{\rho}^{(\alpha)} (x)$. Thus, in this case it is possible to
extend the Laguerre series (\ref{LagSer_ExpoDecaySerCoeffs}), which
converges \emph{exponentially} in the mean for $\vert s \vert < 1$, to
the boundary of its circle of convergence, yielding a \emph{convergent}
Laguerre series with algebraically decaying coefficients that ultimately
have strictly alternating signs.

It is important to note that this approach does not always work. If we
make in the Laguerre series (\ref{calE_LagSer}), which because of
(\ref{calE_PowSer_4}) is nothing but the generating function
(\ref{GLag_Gen_1}), the substitution $t \to -1$, we formally obtain:
\begin{equation}
  \mathrm{e}^{z/2} \; = \; 2^{\alpha+1} \,
  \sum_{n=0}^{\infty} \, (-1)^{n} \, L_{n}^{(\alpha)} (z) \, .
\end{equation}
However, $\mathrm{e}^{z/2}$ does not belong to the Hilbert space
$L^{2}_{z^{\alpha} \mathrm{e}^{-z}} \bigl( [0, \infty) \bigr)$ defined by
(\ref{HilbertL^2_Lag}), and it follows from (\ref{ConvCon_LambdaSer})
that the series (\ref{ConvCond_f_Exp_GenLag_Normalized}) for the norm
$\bigl\Vert \mathcal{\mathrm{e}}^{z/2} \bigr\Vert_{z^{\alpha}
  \mathrm{e}^{-z}, 2}$ does not converge.

It is also not possible to assume $\vert s \vert > 1$ in the Laguerre
series (\ref{LagSer_ExpoDecaySerCoeffs}) for
$\mathcal{G}_{\rho}^{(\alpha)} (s; z)$, because then the series
(\ref{ConvCond_f_Exp_GenLag_Normalized}) for the norm $\bigl\Vert
\mathcal{G}_{\rho}^{(\alpha)} (s; z) \bigr\Vert_{z^{\alpha}
  \mathrm{e}^{-z}, 2}$ diverges. In contrast, the explicit expression
(\ref{Rearr_calG_Exp_GLag_4}) for $\mathcal{G}_{\rho}^{(\alpha)} (s; z)$
remains meaningful for $\vert s \vert > 1$.

If the series coefficients of a power series decay \emph{factorially}, we
can expect that the function represented by this power series is analytic
in the whole complex plane $\mathbb{C}$. The most obvious examples are
the exponential $\exp (z)$, or -- slightly more general -- the confluent
hypergeometric series ${}_{1} F_{1} (a; b; z)$. In contrast, power
series, whose coefficients decay algebraically, converge only in compact
subset of the complex plane. It should be of interest to study how
factorially decaying series coefficients $\lambda_{n}^{(\alpha)}$
influence the analyticity properties of functions represented by a
Laguerre series of the type of (\ref{f_Exp_GLag}).

As an example, let us consider the following Laguerre series:
\begin{equation}
  \label{calF_LagSer}
  \mathcal{F}^{(\alpha)} (s; z) \; = \;
  \sum_{n=0}^{\infty} \, \frac{s^{n}}{(\alpha+1)_{n}} \,
  L_{n}^{(\alpha)} (z) \, .
\end{equation}
Thus, we choose $\lambda_{n}^{(\alpha)} = s^{n}/(\alpha+1)_{n}$ in
(\ref{f_Exp_GLag}). Inserting this into the modified rearranged Laguerre
expansion (\ref{Rearr_f_Exp_GLag_1}) yields:
\begin{align}
  \label{Rearr_calF_Exp_GLag_2}
  \mathcal{F}^{(\alpha)} (s; z) & \; = \;
  \sum_{\nu=0}^{\infty} \, \frac{(-s z)^{\nu}}{(\alpha+1)_{\nu} \, \nu!}
  \, \sum_{\mu=0}^{\infty} \, \frac{s^{\mu}}{\mu!}
  \\
  \label{Rearr_calF_Exp_GLag_3}
  & \; = \; \mathrm{e}^{s} \,
  \sum_{\nu=0}^{\infty} \, \frac{(-s z)^{\nu}}{(\alpha+1)_{\nu} \, \nu!}
\end{align}
The infinite series in (\ref{Rearr_calF_Exp_GLag_2}) can be expressed as
a generalized hypergeometric series ${}_{0} F_{1}$. We then obtain:
\begin{equation}
  \label{Rearr_calF_Exp_GLag_4}
  \mathcal{F}^{(\alpha)} (s; z) \; = \; \mathrm{e}^{s} \,
  {}_{0} F_{1} \bigl( \alpha+1; -s z \bigr) \, .
\end{equation}
However, this is nothing but a known generating function \cite[p.\
242]{Magnus/Oberhettinger/Soni/1966}.

The generalized hypergeometric series ${}_{0} F_{1} \bigl( \alpha+1; -s z
\bigr)$ in (\ref{Rearr_calF_Exp_GLag_4}) obviously converges for all $s,
z \in \mathbb{C}$ as long as $\alpha+1$ is not a negative integer. Since
we always assume $\alpha > -1$, we thus can conclude that
$\mathcal{F}^{(\alpha)} (s; z)$ is an analytic function in the sense of
complex analysis in \emph{both} $s$ and $z$. Moreover, the series
(\ref{MeanConvCond_f}) for the norm $\bigl\Vert \mathcal{F}^{(\alpha)}
(s; z) \bigr\Vert_{z^{\alpha} \mathrm{e}^{-z}, 2}$ converges for
arbitrary $s \in \mathbb{C}$.

\typeout{==> Section: Computational Approaches}
\section{Computational Approaches}
\label{Sec:ComputationalApproaches}

The examples considered in Sections
\ref{Sec:AlgebraicallyDecayingSeriesCoefficients} and
\ref{Sec:ExponentiallyFactoriallyDecayingSeriesCoefficients} show that it
is indeed possible to transform a Laguerre series of the type of
(\ref{f_Exp_GLag}) with the help of (\ref{Rearr_f_Exp_GLag}) to a power
series about $z=0$ of the type of (\ref{PowSer_f}). The essential
requirement is that the inner $\mu$ series in (\ref{Rearr_f_Exp_GLag})
has to converge. This is the case if the coefficients
$\lambda_{n}^{(\alpha)}$ of the Laguerre series decay sufficiently
rapidly, or -- if the coefficients $\lambda_{n}^{(\alpha)}$ only decay
algebraically in magnitude -- ultimately have strictly alternating signs.

In Section \ref{Sec:ExponentiallyFactoriallyDecayingSeriesCoefficients},
the transformation of Laguerre series with exponentially and factorially
decaying series coefficients was discussed. Apart from finite Laguerre
expansions, this is pretty much the best scenario which can occur in this
context. The convergence of the inner $\mu$ series is guaranteed, and
often such a $\mu$ series converges rapidly, making the resulting power
series expansion computationally useful even if it is not possible to
find a convenient closed form expression for the inner $\mu$ series.

As discussed in Section
\ref{Sec:AlgebraicallyDecayingSeriesCoefficients}, the most challenging
and therefore also most interesting problems occur if the coefficients
$\lambda_{n}^{(\alpha)}$ decay algebraically as $n \to \infty$. If
algebraically decaying coefficients $\lambda_{n}^{(\alpha)}$ ultimately
have the same sign, the inner $\mu$ series in (\ref{Rearr_f_Exp_GLag})
diverges, and it is also not possible to sum these series to something
finite. Accordingly, the transformation formula (\ref{Rearr_f_Exp_GLag})
leads to a power series expansion having infinitely many series
coefficients that are infinite in magnitude. This simply means that a
power series in $z$ does not exist because the function under
consideration is not analytic at the origin.

The inner $\mu$ series in (\ref{Rearr_f_Exp_GLag}) also does not converge
if the coefficients $\lambda_{n}^{(\alpha)}$ decay algebraically and
ultimately have strictly alternating signs. But in this case, suitable
summation techniques are capable of associating something finite to the
divergent inner $\mu$ series.

In all examples considered in Section
\ref{Sec:AlgebraicallyDecayingSeriesCoefficients}, it was possible to sum
the divergent inner $\mu$ series with ultimately strictly alternating
signs by means of explicit analytic continuation formulas for
hypergeometric series ${}_{1} F_{0}$ and ${}_{2} F_{1}$, respectively.
For example, the Laguerre series (\ref{LagSer_AltSerCoeffs}) for
$G_{\rho}^{(\alpha)} (z)$ leads according to (\ref{Rearr_G_Exp_GLag_3})
to an inner $\mu$ series that can be expressed as a binomial series
${}_{1} F_{0} (-\rho+\nu; -1)$ that diverges for fixed $\rho \in
\mathbb{R} \setminus \mathbb{N}_{0}$ at least for sufficiently large
values of $\nu \in \mathbb{N}_{0}$. Nevertheless, the summation of this
divergent series is trivial because every hypergeometric series ${}_1 F_0
(a; z)$ with $z \ne 1$ possesses according to (\ref{BinomSer}) a very
simple closed form expression which accomplishes the necessary analytic
continuation.

Similarly, the Laguerre series (\ref{LagSer_2H1_a_b_c_alpha}) for
$\prescript{}{2}{\mathcal{H}}_{1}^{(\alpha)} (a, b; c; z)$ leads
according to (\ref{Rearr_2H1_Exp_GLag_3}) to an inner $\mu$ series that
can be expressed as a hypergeometric series ${}_{2} F_{1} ( a+\nu, b+\nu;
c+\nu; -1)$ that also diverges for sufficiently large $\nu \in
\mathbb{N}_{0}$. But again, it is almost trivially simple to find
explicit analytic continuation formulas that replace this divergent
series by convergent expansions. The linear transformations (\ref{LTr_1})
and (\ref{LTr_2}) do the job.

Unfortunately, the situation is not nearly as nice if the series
coefficients $\lambda_{n}^{(\alpha)}$ have a more complicated structure.
Let us for example consider the following Laguerre series with
coefficients $\lambda_{n}^{(\alpha)}$ that are ratios of $p+1$ numerator
and $p+1$ denominator Pochhammer symbols:
\begin{align}
  \label{LagSer_(p+1)_H_p_alpha}
  & \prescript{}{p+1}{\mathcal{H}}_{p}^{(\alpha)}
  \bigl(a_{1}, \dots, a_{p+1}; b_{1}, \dots, b_{p}; z \bigr)
  \notag \\
  & \qquad
  \; = \; \sum_{n=0}^{\infty} \, (-1)^{n} \,
  \frac{(a_{1})_{n} \dots (a_{p+1})_{n}}{(b_{1})_{n} \dots (b_{p})_{n}}
  \, \frac{L_{n}^{(\alpha)} (z)}{(\alpha+1)_{n}} \, ,
  \qquad p \in \mathbb{N}_{0} \, .
\end{align}
Obviously, this Laguerre series generalizes the Laguerre series
(\ref{LagSer_2H1_a_b_c_alpha}) for
$\prescript{}{2}{\mathcal{H}}_{1}^{(\alpha)} (a, b; c; z)$.

It is easy to show that the application of the transformation formula
(\ref{Rearr_f_Exp_GLag}) to the Laguerre series
(\ref{LagSer_(p+1)_H_p_alpha}) leads to inner $\mu$ series that can be
expressed as generalized hypergeometric series
\begin{equation}
  \label{Div_(p+1)_F_p}
  {}_{p+1} F_{p} \left( \genfrac{}{}{0pt}{}
  {a_{1}+\nu, \dots, a_{p+1}+\nu}
  {b_{1}+\nu, \dots, b_{p}+\nu}; z \right)
  \; = \; \sum_{n=0}^{\infty} \,
  \frac{(a_{1}+\nu)_{n} \cdots (a_{p+1}+\nu)_{n}}
  {(b_{1}+\nu)_{n} \cdots (b_{p}+\nu)_{n}} \, \frac{z^{n}}{n!}
\end{equation}
with argument $z=-1$. Just like the Gaussian hypergeometric series
${}_{2} F_{1}$, a nonterminating generalized hypergeometric series
${}_{p+1} F_{p}$ converges only in the interior of the unit circle.
Moreover, (\ref{AsyGammaRatio}) implies
\begin{equation}
  \frac{\Gamma (a_{1}+\nu+n) \, \Gamma (a_{p+1}+\nu+n)}
  {\Gamma (b_{1}+\nu+n) \, \Gamma (b_{p}+\nu+n) \, n!} \; \sim \;
  n^{a_{1} + \dots + a_{p+1} + \nu - b_{1} + \dots - b_{p}} \, ,
  \qquad n \to \infty \, .
\end{equation}
This asymptotic estimate shows that the generalized hypergeometric series
(\ref{Div_(p+1)_F_p}) with argument $z=-1$ diverge for sufficiently large
values of $\nu$ and have to be summed. Unfortunately, the theory of the
generalized hypergeometric series ${}_{p+1} F_{p}$ with $p \ge 2$ is not
nearly as highly developed as the theory of the Gaussian hypergeometric
series ${}_{2} F_{1}$.  Many transformation and/or analytic continuation
formulas for a generalized hypergeometric series ${}_{p+1} F_{p}$ with $p
\ge 2$ are either not known or at least much more complicated than the
corresponding expressions for a Gaussian hypergeometric series ${}_{2}
F_{1}$ (compare also \cite{Buehring/1988}).  We can safely assume that
our chances of finding convenient explicit analytic continuation
formulas, that can accomplish the summation of the divergent generalized
hypergeometric series (\ref{Div_(p+1)_F_p}), are rapidly approaching zero
with increasing $p$.

Thus, the rearrangement of Laguerre series via (\ref{Rearr_f_Exp_GLag})
can only produce closed form expressions if algebraically decaying and
ultimately strictly alternating series coefficients
$\lambda_{n}^{(\alpha)}$ have an exceptionally simple structure, as for
example the ones in the Laguerre series (\ref{LagSer_AltSerCoeffs}) and
(\ref{LagSer_2H1_a_b_c_alpha}). In the case of more complicated
coefficients $\lambda_{n}^{(\alpha)}$, the necessary summations of the
divergent inner $\mu$ series will be too difficult to produce closed form
expressions of manageable complexity.

These and related considerations may tempt a skeptical reader to argue
with some justification that all examples considered in Section
\ref{Sec:AlgebraicallyDecayingSeriesCoefficients} are actually fairly
simple. Therefore, it is by no means obvious whether the approach of this
article is capable of producing anything new beyond a rederivation of
known and comparatively simple formulas.

While the first conclusion is certainly correct, the second is in my
opinion overly pessimistic. The reason is that summation via analytic
continuation formulas is not the only possibility: We can also try to use
\emph{numerical} summation techniques. Of course, knowing only numerical
approximations to the leading coefficients of a power series is not
nearly as nice as knowing explicit and possibly even simple expressions
for the series coefficients, but it is certainly better than nothing. In
practical applications, this limited information may suffice.

The use of numerical summation techniques is by no means a new idea.
Already in 1882, H\"{o}lder \cite{Hoelder/1882} proposed -- inspired by
an article by Frobenius \cite{Frobenius/1880} -- to determine the value
of a power series on the boundary of its circle of convergence with the
help of numerical summation processes based on weighted arithmetic means.

As discussed in Appendix \ref{App:SequenceTransformations}, the numerical
processes proposed by H\"{o}lder \cite{Hoelder/1882} and others
subsequently evolved to a sophisticated mathematical theory of so-called
regular matrix transformations that have many advantageous theoretical
features and that can be used for the summation of divergent alternating
series possessing similar divergence properties as the generalized
hypergeometric series (\ref{Div_(p+1)_F_p}).

Unfortunately, regular matrix transformations are in general at best
moderately powerful. Therefore, I prefer to use instead nonlinear
sequence transformations, whose theoretical properties are not yet
completely understood, but which are according to experience often much
better suited to achieve highly accurate summation results.

It may be interesting to note that summation techniques such as conformal
mappings, reexpansions, and Pad\'{e} approximants were discussed in an
article by Skorokhodov \cite{Skorokhodov/2003} on the analytic
continuation of a divergent generalized hypergeometric series ${}_{p+1}
F_{p}$ with argument $\vert z \vert \ge 1$. Skorokhodov completely
ignored the theoretically much simpler, but also less powerful regular
matrix transformations. Nevertheless, I do not think that the summation
techniques considered by Skorokhodov give best results in the case of the
divergent, but summable inner $\mu$ series occurring in this article. In
my opinion, (much) better results can be obtained with the help of those
nonlinear sequence transformations that are reviewed in Appendix
\ref{App:SequenceTransformations}.

In this article, I apply as numerical summation techniques Wynn's epsilon
algorithm (\ref{eps_al}), and the two Levin-type transformations
(\ref{dLevTr}) and (\ref{dWenTr}). Wynn's epsilon algorithm produces
Pad\'{e} approximants according to (\ref{Eps_Pade}) if the input data are
the partial sums $f_{n} (z) = \sum_{k=0}^{n} \gamma_{k} z^{k}$ of a power
series. The two Levin-type transformations (\ref{dLevTr}) and
(\ref{dWenTr}) are based on the remainder estimate (\ref{dRemEst}) which
corresponds to the first term neglected in the partial sum. Both are
known to be highly effective in the case of both convergent and divergent
alternating series.

The discussion in Section
\ref{Sec:AlgebraicallyDecayingSeriesCoefficients} should provide
convincing evidence that for $z=-1$ and in particular for large values of
$\nu \in \mathbb{N}_{0}$ it is not a particularly good idea to use the
hypergeometric series ${}_{2} F_{1} (a+\nu, b+\nu; c+\nu; z)$ for the
evaluation of the hypergeometric function it represents. With the help of
analytic continuation formulas like (\ref{LTr_1}) and (\ref{LTr_2}) or
also (\ref{LTr_0}), computationally much more convenient hypergeometric
series can be derived. However, these simplifications are only possible
in the case of a Gaussian hypergeometric series ${}_{2} F_{1}$, but not
necessarily in the case of more complicated generalized hypergeometric
series ${}_{p+1} F_{p}$, let alone in the case of divergent, but summable
inner $\mu$ series that result from purely numerical Laguerre series
coefficients $\lambda_{n}^{(\alpha)}$. Thus, the hypergeometric series
${}_{2} F_{1} (a+\nu, b+\nu; c+\nu; -1)$ serves as a model problem for
other divergent alternating series whose terms also increase in magnitude
like a fixed power of the index. It should be interesting to see how much
can be accomplished by employing powerful nonlinear sequence
transformations.

In Tables \ref{Tab_6_1} and \ref{Tab_6_2}, the nonlinear transformations
mentioned above are applied to the partial sums
\begin{equation}
  \label{ParSum_2F1}
  s_{n} \; = \; s_{n} (a, b, c, \nu) \; = \; \sum_{k=0}^{n} \,
  (-1)^{k} \, \frac{(a+\nu)_{k} (b+\nu)_{k}}{(c+\nu)_{k} k!}
\end{equation}
of the Gaussian hypergeometric series ${}_{2} F_{1} (a+\nu, b+\nu; c+\nu;
-1)$ with $a=3/2$, $b=7/3$, and $c=21/4$. In Table \ref{Tab_6_1} we have
$\nu=0$, and in Table \ref{Tab_6_2} we have $\nu=10$. All transformation
results of this article were obtained using the floating point
arithmetics of Maple 11, and the ``exact'' result in Tables \ref{Tab_6_1}
and \ref{Tab_6_2} was produced by the Maple procedure \texttt{hypergeom}
which computes generalized hypergeometric series.

The entries in the third column of Table \ref{Tab_6_1} are chosen
according (\ref{EpsAl_ApprLim_gen}), and the entries in columns 4 and 5
according to (\ref{LevinTypeApprLim}). Since Wynn's epsilon algorithm
produces according to (\ref{Eps_Pade}) Pad\'{e} approximants if the input
data are the partial sums of a power series, the entries in the third
column can be identified with the staircase sequence $\bigl[ \Ent
{(n+1)/2}/\Ent {n/2} \bigr]$ of Pad\'{e} approximants to the
hypergeometric series ${}_2 F_{1} (a+\nu, b+\nu; c+\nu; z)$ with $z=-1$.
Here, $\Ent {x}$ stands for the integral part of $x$, which is the
largest integer $\nu$ satisfying $\nu \le x$.

\begin{table}[t]
  \caption{Application of Wynn's epsilon algorithm and the
    Levin-type transformations $d_{k}^{(n)} (\beta, s_n)$ and
    $\delta_{k}^{(n)} (\beta, s_n)$ with $\beta=1$ to the partial sums
    of the slowly convergent hypergeometric series
    ${}_2 F_{1} (a+\nu, b+\nu; c+\nu; -1)$ with $a = 3/2$, $b = 7/3$,
    $c = 21/4$, and $\nu = 0$.}
  \label{Tab_6_1}
  \begin{tabular*}{\textwidth}{lrrrr}%
    \\
    \hline \hline %
    $n$%
    & \multicolumn{1}{c}{$s_{n} (z)$}%
    & \multicolumn{1}{c}{$\epsilon_{2 \Ent {n/2}}^{(n - 2 \Ent {n/2})}$}%
    & \multicolumn{1}{c}{$d_{n}^{(0)} \bigl(1, s_0 (z) \bigr)$}%
    & \multicolumn{1}{c}{${\delta}_{n}^{(0)} \bigl(1, s_0 (z) \bigr)$}%
    \rule[-4pt]{0pt}{20pt} \\
    & \multicolumn{1}{c}{Eq.\ (\protect\ref{ParSum_2F1})}%
    & \multicolumn{1}{c}{Eq.\ (\protect\ref{eps_al})}%
    & \multicolumn{1}{c}{Eq.\ (\protect\ref{dLevTr})}%
    & \multicolumn{1}{c}{Eq.\ (\protect\ref{dWenTr})}%
    \rule[-6pt]{0pt}{12pt}
    \\
    \hline%
    0  & 1.00000000 & 1.000000000000000 & 1.000000000000000
    & 1.000000000000000 \rule[-1pt]{0pt}{12pt} \\
    1  & 0.33333333 & 0.333333333333333 & 0.600000000000000
    & 0.600000000000000 \rule[-1pt]{0pt}{12pt} \\
    2  & 0.77777778 & 0.600000000000000 & 0.596079232182969
    & 0.596079232182969 \rule[-1pt]{0pt}{12pt} \\
    3  & 0.46785866 & 0.595184349134688 & 0.597361776678518
    & 0.597283629053888 \rule[-1pt]{0pt}{12pt} \\
    4  & 0.69325438 & 0.597114931459132 & 0.597128156362257
    & 0.597151427941970 \rule[-1pt]{0pt}{12pt} \\
    5  & 0.52349689 & 0.597142098567676 & 0.597159499152400
    & 0.597156391214500 \rule[-1pt]{0pt}{12pt} \\
    6  & 0.65507045 & 0.598028964909470 & 0.597156082629105
    & 0.597156376752544 \rule[-1pt]{0pt}{12pt} \\
    7  & 0.55064699 & 0.597156986187236 & 0.597156397164508
    & 0.597156374043337 \rule[-1pt]{0pt}{12pt} \\
    8  & 0.63518026 & 0.597156454880736 & 0.597156372412882
    & 0.597156373980610 \rule[-1pt]{0pt}{12pt} \\
    9  & 0.56559242 & 0.597156266466529 & 0.597156374069016
    & 0.597156373980877 \rule[-1pt]{0pt}{12pt} \\
    10 & 0.62370437 & 0.597156373786525 & 0.597156373977124
    & 0.597156373980968 \rule[-1pt]{0pt}{12pt} \\
    11 & 0.57457047 & 0.597156373530068 & 0.597156373981079
    & 0.597156373980973 \rule[-1pt]{0pt}{12pt} \\
    \hline
    exact &         & 0.597156373980973 & 0.597156373980973
    & 0.597156373980973 \rule[-1pt]{0pt}{12pt} \\
    \hline \hline %
  \end{tabular*}
\end{table}

The asymptotic estimate (\ref{AsyGammaRatio}) shows that the terms
$(-1)^{n} (a)_{n} (b)_{n}/[{c}_{n} n!]$ of the hypergeometric series in
Table \ref{Tab_6_1} decay in magnitude like $n^{-29/12}$ as $n \to
\infty$. Therefore, this hypergeometric series converges, albeit quite
slowly. This assessment is confirmed by the data in the second column of
Table \ref{Tab_6_1}. The transformation results in Table \ref{Tab_6_1}
indicate that all three transformations are able to accelerate the
convergence of the hypergeometric series ${}_{2} F_{1}$ effectively. The
least effective, but still very powerful accelerator is Wynn's epsilon
algorithm, and the most effective transformation is the Levin-type delta
transformation which is -- as documented by numerous references mentioned
in Appendix \ref{App:SequenceTransformations} -- known to be highly
effective in the case of both convergent and divergent alternating
series.

\begin{table}[t]
  \caption{Summation of the hypergeometric series
    ${}_2 F_{1} (a+\nu, b+\nu; c+\nu; -1)$ with $a = 3/2$, $b = 7/3$,
    $c = 21/4$, and $\nu = 10$ with the help of Wynn's epsilon
    algorithm and the Levin-type transformation
    $\delta_{k}^{(n)} (\beta, s_n)$ with $\beta=1$.}
  \label{Tab_6_2}
  \begin{tabular*}{\textwidth}{lrrr}%
    \\
    \hline \hline %
    $n$%
    & \multicolumn{1}{c}{$s_{n} (z)$}%
    & \multicolumn{1}{c}{$\epsilon_{2 \Ent {n/2}}^{(n - 2 \Ent {n/2})}$}%
    & \multicolumn{1}{c}{${\delta}_{n}^{(0)} \bigl(1, s_0 (z) \bigr)$}%
    \rule[-4pt]{0pt}{20pt} \\
    & \multicolumn{1}{c}{Eq.\ (\protect\ref{ParSum_2F1})}%
    & \multicolumn{1}{c}{Eq.\ (\protect\ref{eps_al})}%
    & \multicolumn{1}{c}{Eq.\ (\protect\ref{dWenTr})}%
    \rule[-6pt]{0pt}{12pt}
    \\
    \hline%
    0  & $ 0.10000 \times 10^{1}$ & $ 0.100000000000000 \times 10^{+1}$
    & $ 0.100000000000000 \times 10^{+1}$ \\
    1  & $-0.83005 \times 10^{1}$ & $-0.830054644808743 \times 10^{+1}$
    & $-0.517662391110501 \times 10^{+0}$ \\
    2  & $ 0.39395 \times 10^{2}$ & $-0.517662391110501 \times 10^{+0}$
    & $ 0.355800262535272 \times 10^{+0}$ \\
    3  & $-0.13894 \times 10^{3}$ & $ 0.176355932965950 \times 10^{+1}$
    & $-0.979095090620413 \times 10^{-1}$ \\
    4  & $ 0.40421 \times 10^{3}$ & $ 0.150567485119558 \times 10^{+0}$
    & $ 0.156275192732603 \times 10^{-1}$ \\
    5  & $-0.10245 \times 10^{4}$ & $-0.292708508243714 \times 10^{+0}$
    & $ 0.539203120343906 \times 10^{-3}$ \\
    6  & $ 0.23385 \times 10^{4}$ & $-0.235266825737199 \times 10^{-1}$
    & $ 0.137225149988457 \times 10^{-2}$ \\
    7  & $-0.49149 \times 10^{4}$ & $ 0.351057409892771 \times 10^{-1}$
    & $ 0.145137993614918 \times 10^{-2}$ \\
    8  & $ 0.96599 \times 10^{4}$ & $ 0.394766824723114 \times 10^{-2}$
    & $ 0.142550654190640 \times 10^{-2}$ \\
    9  & $-0.17957 \times 10^{5}$ & $-0.113107334634727 \times 10^{-2}$
    & $ 0.142932832720914 \times 10^{-2}$ \\
    10 & $ 0.31849 \times 10^{5}$ & $ 0.127361316374322 \times 10^{-2}$
    & $ 0.142892054773978 \times 10^{-2}$ \\
    11 & $-0.54255 \times 10^{5}$ & $ 0.155543251021275 \times 10^{-2}$
    & $ 0.142895407283313 \times 10^{-2}$ \\
    12 & $ 0.89251 \times 10^{5}$ & $ 0.143471035609870 \times 10^{-2}$
    & $ 0.142895196573040 \times 10^{-2}$ \\
    13 & $-0.14240 \times 10^{6}$ & $ 0.142504106397527 \times 10^{-2}$
    & $ 0.142895206303027 \times 10^{-2}$ \\
    14 & $ 0.22113 \times 10^{6}$ & $ 0.142883205557579 \times 10^{-2}$
    & $ 0.142895205998921 \times 10^{-2}$ \\
    15 & $-0.33524 \times 10^{6}$ & $ 0.142902170361815 \times 10^{-2}$
    & $ 0.142895206004165 \times 10^{-2}$ \\
    16 & $ 0.49740 \times 10^{6}$ & $ 0.142895322013589 \times 10^{-2}$
    & $ 0.142895206004153 \times 10^{-2}$ \\
    17 & $-0.72381 \times 10^{6}$ & $ 0.142895148436317 \times 10^{-2}$
    & $ 0.142895206004152 \times 10^{-2}$ \\
    \hline
    exact &                      & $ 0.142895206004152 \times 10^{-2}$
    & $ 0.142895206004152 \times 10^{-2}$ \rule[-1pt]{0pt}{12pt} \\
    \hline \hline %
  \end{tabular*}
\end{table}

In Table \ref{Tab_6_2}, only Wynn's epsilon algorithm and the delta
transformation are displayed. The data in the second column show that the
partial sums (\ref{ParSum_2F1}) with $\nu=10$ diverge rapidly.
Nevertheless, it is possible to obtain highly accurate summation results.
As in Table \ref{Tab_6_1}, the delta transformation was clearly more
effective than Wynn's epsilon algorithm.

It is in my opinion remarkable that Levin's transformation
(\ref{dLevTr}), whose action on the partial sums (\ref{ParSum_2F1}) is
not displayed in Table \ref{Tab_6_2}, turned out to be only roughly as
effective as Wynn's epsilon algorithm. Under the same conditions as in
Table \ref{Tab_6_2}, I obtained the following summation results:
\begin{align}
  d_{16}^{(0)} \bigl(1, s_0 (z) \bigr) & \; = \;
  0.142894441246992 \times 10^{-2} \, ,
  \\
  d_{17}^{(0)} \bigl(1, s_0 (z) \bigr) & \; = \;
  0.142895355463039 \times 10^{-2} \, .
\end{align}
In view of the in general very good reputation of Levin's sequence
transformation, this comparatively weak performance is somewhat puzzling.

In the case of more complicated generalized hypergeometric series
${}_{p+1} F_{p}$, the same general pattern was observed. If the same
nonlinear sequence transformations as in Tables \ref{Tab_6_1} and
\ref{Tab_6_2} are applied to the partial sums
\begin{equation}
  \label{ParSum_3F2}
  s_{n} \; = \; s_{n} (a, b, c, d, e, \nu) \; = \; \sum_{k=0}^{n} \,
  (-1)^{k} \, \frac{(a+\nu)_{k} (b+\nu)_{k} (c+\nu)_{k}}
  {(d+\nu)_{k} (e+\nu)_{k} k!}
\end{equation}
of the generalized hypergeometric series ${}_{3} F_{2} (a+\nu, b+\nu,
c+\nu; d+\nu, e+\nu; -1)$ with $a = 3/2$, $b = 7/3$, $c = 11/5$, $d =
22/7$, and $e = 32/11$, we obtain for the convergent hypergeometric
series with $\nu=0$
\begin{align}
  \epsilon_{12}^{(0)}
  & \; = \; 0.536266961325332 \, ,
  \\
  d_{12}^{(0)} \bigl(1, s_0 (z) \bigr)
  & \; = \; 0.536266961240988 \, ,
  \\
  {\delta}_{12}^{(0)} \bigl(1, s_0 (z) \bigr)
  & \; = \; 0.536266961240986 \, ,
  \\
  \texttt{hypergeom} & \; = \; 0.536266961240986 \, ,
\end{align}
and for the divergent hypergeometric series with $\nu=10$
\begin{align}
  \epsilon_{16}^{(1)}
  & \; = \; 0.816454762672306 \times 10^{-3} \, ,
  \\
  d_{17}^{(0)} \bigl(1, s_0 (z) \bigr)
  & \; = \; 0.816459448502108 \times 10^{-3} \, ,
  \\
  {\delta}_{17}^{(0)} \bigl(1, s_0 (z) \bigr)
  & \; = \; 0.816458731770118 \times 10^{-3} \, ,
  \\
  \texttt{hypergeom} & \; = \; 0.816458731770118 \times 10^{-3} \, .
\end{align}

In the case of the partial sums
\begin{equation}
  \label{ParSum_4F3}
  s_{n} \; = \; s_{n} (a, b, c, d, e, f, g \nu) \; = \;
  \sum_{k=0}^{n} \,
  (-1)^{k} \, \frac{(a+\nu)_{k} (b+\nu)_{k} (c+\nu)_{k} (d+\nu)_{k}}
  {(e+\nu)_{k} (f+\nu)_{k} (g+\nu)_{k} k!}
\end{equation}
of the generalized hypergeometric series ${}_{4} F_{3} (a+\nu, b+\nu,
c+\nu, d+\nu; e+\nu, f+\nu, g+\nu; -1)$ with $a = 3/2$, $b = 7/3$, $c =
11/5$, $d = 16/17$, $e = 18/19$, $f = 22/7$, and $g = 32/11$, we obtain
for the convergent hypergeometric series with $\nu=0$
\begin{align}
  \epsilon_{12}^{(0)}
  & \; = \; 0.538509188429164 \, ,
  \\
  d_{12}^{(0)} \bigl(1, s_0 (z) \bigr)
  & \; = \; 0.538509188330837 \, ,
  \\
  {\delta}_{12}^{(0)} \bigl(1, s_0 (z) \bigr)
  & \; = \; 0.538509188330835 \, ,
  \\
  \texttt{hypergeom} & \; = \; 0.538509188330835 \, ,
\end{align}
and for the divergent hypergeometric series with $\nu=10$
\begin{align}
  \epsilon_{16}^{(1)}
  & \; = \; 0.819691553980977 \times 10^{-3} \, ,
  \\
  d_{17}^{(0)} \bigl(1, s_0 (z) \bigr)
  & \; = \; 0.819696193819277 \times 10^{-3} \, ,
  \\
  {\delta}_{17}^{(0)} \bigl(1, s_0 (z) \bigr)
  & \; = \; 0.819695479036364 \times 10^{-3} \, ,
  \\
  \texttt{hypergeom} & \; = \; 0.819695479036364 \times 10^{-3} \, .
\end{align}

It was emphasized in Section
\ref{Sec:ExponentiallyFactoriallyDecayingSeriesCoefficients} that the
inner $\mu$ series in (\ref{Rearr_f_Exp_GLag}) converge if the
coefficients $\lambda_{n}^{(\alpha)}$ of a Laguerre series decay
exponentially as $n \to \infty$. Accordingly, summation techniques
are not needed in this case. Nevertheless, nonlinear sequence
transformations can be extremely useful even in the case of exponentially
decaying Laguerre series coefficients. Let us for example assume that the
coefficients $\lambda_{n}^{(\alpha)}$ of a Laguerre series satisfy
\begin{equation}
  \lambda_{n}^{(\alpha)} \; = \; t^{n} \, \ell_{n}^{(n)}
\end{equation}
and that it is not possible to derive a closed form expression for the
inner $\mu$ series in (\ref{Rearr_f_Exp_GLag}). Thus, the the inner $\mu$
series have to be evaluated numerically.

Let us now also assume that the coefficients $\ell_{n}^{(n)}$
\emph{increase} like a fixed power of $n$ as $n \to \infty$. In spite of
this unfavorable behavior, the inner $\mu$ series converge as long as
$\vert t \vert < 1$.  However, convergence can become prohibitively slow
if $\vert t \vert$ is only slightly smaller than one. The convergence
problems are particularly severe if the terms in inner $\mu$ series
ultimately have the same sign because then the $\mu$ series do not
converge for $t=1$ and are also not summable. But again, nonlinear
sequence transformations can be extremely useful to speed up the
convergence of such a monotone series with (very) slowly decaying terms
(see for example \cite{Jentschura/Mohr/Soff/Weniger/1999} or
\cite[Section
2.2.6]{Caliceti/Meyer-Hermann/Ribeca/Surzhykov/Jentschura/2007} and
references therein).

As documented by the recent books by Cuyt, Brevik Petersen, Verdonk,
Waadeland, and Jones
\cite{Cuyt/BrevikPetersen/Verdonk/Waadeland/Jones/2008} and by Gil,
Segura, and Temme \cite{Gil/Segura/Temme/2007}, or by a review by Temme
\cite{Temme/2007}, there is currently a lot of work being done on the
efficient and reliable evaluation of special functions. As I had shown in
several articles \cite{Jentschura/Gies/Valluri/Lamm/Weniger/2002,%
  Jentschura/Mohr/Soff/Weniger/1999,Weniger/1989,Weniger/1990,%
  Weniger/1994a,Weniger/1994b,Weniger/1996d,Weniger/2001,Weniger/2003,%
  Weniger/Cizek/1990,Weniger/Steinborn/1989a}, nonlinear sequence
transformations can be extremely useful in this respect.

\typeout{==> Section: Guseinov's Rearranged One-Range Addition Theorems}
\section{Guseinov's Rearranged One-Range Addition Theorems}
\label{Sec:GuseinovsRearrangedOne-RangeAdditionTheorems}

The analyticity of Laguerre series is a problem of classical analysis,
but I became interested in this mathematical topic because of some open
questions in molecular electronic structure theory. During the work for
my forthcoming review on addition theorems \cite{Weniger/2008a*}, I came
across some articles by Guseinov
\cite{Guseinov/1980a,Guseinov/2001a,Guseinov/2002c} who had constructed
one-range addition theorems for Slater-type functions
\cite{Slater/1930,Slater/1932}. In unnormalized form, Slater-type
functions, which play a major role as basis functions in atomic and
molecular electronic structure calculations, are defined as follows:
\begin{equation}
  \label{Def_STF}
\chi_{N, L}^{M} (\beta, \mathbf{r}) \; = \;
(\beta r)^{N-L-1} \, \mathrm{e}^{- \beta r} \,
\mathcal{Y}_{L}^{M} (\beta \bm{r}) \, .
\end{equation}
Here, $\bm{r} \in \mathbb{R}^{3}$, $\mathcal{Y}_{L}^{M} (\beta \bm{r}) =
(\beta r)^{L} Y_{L}^{M} (\theta, \phi)$ is a regular solid harmonic and
$Y_{L}^{M} (\theta, \phi)$ is a (surface) spherical harmonic, $\beta > 0$
is a scaling parameter, $N \in \mathbb{R} \setminus \mathbb{N}$ is a kind
of generalized principal quantum number which is often, but not always a
positive integer $\geq L+1$, and $L$ and $M$ are the usual (orbital)
angular momentum quantum numbers.

Let us assume that $\{ \varphi_{n, \ell}^{m} (\bm{r}) \}_{n, \ell, m}$ is a
\emph{complete} and \emph{orthonormal} function set in the Hilbert space
\begin{equation}
  \label{HilbertL^2}
  L^{2} (\mathbb{R}^3) \; = \; \Bigl\{ f \colon \mathbb{R}^3 \to
  \mathbb{C} \Bigm\vert \, \int \, \vert f (\bm{r}) \vert^2 \,
  \mathrm{d}^3 \bm{r} < \infty \Bigr\}
\end{equation}
of functions that are square integrable with respect to an integration
over the whole $\mathbb{R}^{3}$. Since any $f \in L^{2} (\mathbb{R}^{3})$
can be expanded in terms of the complete and orthonormal functions $\{
\varphi_{n, \ell}^{m} (\bm{r}) \}_{n, \ell, m}$, a one-range addition
theorem for $f (\bm{r} \pm \bm{r}')$ can be formulated as follows:
\begin{subequations}
  \label{OneRangeAddTheor}
  \begin{align}
     \label{OneRangeAddTheor_a}
    f (\bm{r} \pm \bm{r}') & \; = \; \sum_{n \ell m} \, C_{n, \ell}^{m}
    (f; \pm \bm{r}') \, \varphi_{n, \ell}^{m} (\bm{r}) \, ,
    \\
     \label{OneRangeAddTheor_b}
    C_{n, \ell}^{m} (f; \pm \bm{r}') & \; = \; \int \, \bigl[
    \varphi_{n, \ell}^{m} (\bm{r}) \bigr]^{*} \, f (\bm{r} \pm \bm{r}')
    \, \mathrm{d}^3 \bm{r} \, .
  \end{align}
\end{subequations}
The expansion (\ref{OneRangeAddTheor}), which converges in the mean with
respect to the norm of the Hilbert space $L^2 (\mathbb{R}^3)$, is a
one-range addition theorem since the variables ${\bm{r}}$ and ${\bm{r}}'$
are completely separated: The dependence on $\bm{r}$ is entirely
contained in the functions $\varphi_{n, \ell}^{m} (\bm{r})$, whereas
$\bm{r}'$ occurs only in the expansion coefficients $C_{n, \ell}^{m} (f;
\pm \bm{r}')$ which are overlap or convolution-type integrals.

One-range addition theorems of the kind of (\ref{OneRangeAddTheor}) were
constructed by Filter and Steinborn \cite[Eqs.\ (5.11) and
(5.12)]{Filter/Steinborn/1980} and later applied by Kranz and Steinborn
\cite{Kranz/Steinborn/1982} and by Trivedi and Steinborn
\cite{Trivedi/Steinborn/1982}. An alternative derivation of these
addition theorems based on Fourier transformation combined with weakly
convergent expansions of the plane wave $\exp (\pm \mathrm{i} \bm{p}
\cdot \bm{r})$ with $\bm{p}, \bm{r} \in \mathbb{R}^{3}$ was presented in
\cite[Section VII]{Weniger/1985}.

As discussed in \cite[Section 3]{Weniger/2007b}), it is also possible to
formulate one-range addition theorems that converge with respect to the
norm of a \emph{weighted} Hilbert space
\begin{equation}
  \label{HilbertL_w^2}
  L_{w}^{2} (\mathbb{R}^3) \; = \; \Bigl\{ f \colon \mathbb{R}^3 \to
  \mathbb{C} \Bigm\vert \, \int \, w (\bm{r}) \, \vert f (\bm{r}) \vert^2
  \, \mathrm{d}^3 \bm{r} < \infty \Bigr\} \, ,
\end{equation}
where $w (\bm{r}) \ne 1$ is a suitable \emph{positive} weight function.
If we assume that $f \in L_{w}^{2} (\mathbb{R}^3)$ and
that the functions $\{ \psi_{n, \ell}^{m} (\bm{r}) \}_{n, \ell, m}$ are
complete and orthonormal in $L_{w}^{2} (\mathbb{R}^3)$, then we obtain
the following one-range addition theorem \cite[Eq.\
(3.6)]{Weniger/2007b}):
\begin{subequations}
  \label{OneRangeAddTheor_w}
  \begin{align}
    \label{OneRangeAddTheor_w_a}
    f (\bm{r} \pm \bm{r}') & \; = \; \sum_{n \ell m} \, \mathbf{C}_{n,
      \ell}^{m} (f, w; \pm \bm{r}') \, \psi_{n, \ell}^{m} (\bm{r}) \, ,
    \\
    \label{OneRangeAddTheor_w_b}
    \mathbf{C}_{n, \ell}^{m} (f, w; \pm \bm{r}') & \; = \; \int \, \bigl[
    \psi_{n, \ell}^{m} (\bm{r}) \bigr]^{*} \, w (\bm{r}) \, f (\bm{r} \pm
    \bm{r}') \, \mathrm{d}^3 \bm{r} \, .
  \end{align}
\end{subequations}

A one-range addition theorem for a function $f \colon \mathbb{R}^{3} \to
\mathbb{C}$ is a mapping $\mathbb{R}^3 \times \mathbb{R}^3 \to
\mathbb{C}$. Compared to the better known two-range addition theorems
like the so-called Laplace expansion of the Coulomb or Newton potential
$1/r$, which possesses a characteristic two-range form (see for example
\cite[Eq.\ (1.2)]{Weniger/2007b}), one-range addition theorems have the
highly advantageous feature that they provide a \emph{unique} infinite
series representation of $f (\bm{r} \pm \bm{r}')$ with \emph{separated}
variables $\bm{r}$ and $\bm{r}'$ that is valid for the \emph{whole}
argument set $\mathbb{R}^3 \times \mathbb{R}^3$. Further properties of
addition theorems in general and of one-range addition theorems in
particular will be discussed in my forthcoming review
\cite{Weniger/2008a*}.

In his one-range addition theorems, Guseinov used as a complete and
orthonormal function set the following functions \cite[Eq.\
(1)]{Guseinov/2002b}, which -- if the mathematical notation for the
generalized Laguerre polynomials is used -- can be expressed as follows
\cite[Eq.\ (4.16)]{Weniger/2007b}:
\begin{align}
  \label{Def_Psi_Guseinov}
  \prescript{}{k}{\Psi}_{n, \ell}^{m} (\gamma, \bm{r}) & \; = \; \left[
  \frac{(2\gamma)^{k+3} (n-\ell-1)!}{(n+\ell+k+1)!} \right]^{1/2} \,
  \mathrm{e}^{-\gamma r} \, L_{n-\ell-1}^{(2\ell+k+2)} (2 \gamma r) \,
  \mathcal{Y}_{\ell}^{m} (2 \gamma \bm{r}) \, ,
  \notag \\
  & \qquad n \in \mathbb{N} \, , \quad k = -1, 0, 1, 2, \dots \, ,
   \quad \gamma > 0 \, .
\end{align}
As discussed in the text following \cite[Eq.\ (4.20)]{Weniger/2007b},
Guseinov's functions (\ref{Def_Psi_Guseinov}) can -- depending on the
value of $k = -1, 0, 1, 2, \dots$ -- reproduce several other physically
relevant complete and orthonormal function sets.

Guseinov's functions are orthonormal with respect to the weight function
$w (\bm{r}) = r^{k}$ (compare also \cite[Eq.\ (4)]{Guseinov/2002c}):
\begin{equation}
  \label{Psi_Guseinov_OrthoNor}
  \int \, \bigl[ \prescript{}{k}{\Psi}_{n, \ell}^{m} (\gamma, \bm{r})
  \bigr]^{*} \, r^k \, \prescript{}{k}{\Psi}_{n', \ell'}^{m'} (\gamma,
  \bm{r}) \, \mathrm{d}^3 \bm{r} \; = \;
  \delta_{n n'} \, \delta_{\ell \ell'} \, \delta_{m m'} \, .
\end{equation}
Accordingly, Guseinov's functions are complete and orthonormal in the
weighted Hilbert space
\begin{equation}
  \label{HilbertL_r^k^2}
  L_{r^k}^{2} (\mathbb{R}^3) \; = \; \Bigl\{ f \colon \mathbb{R}^3 \to
  \mathbb{C} \Bigm\vert \, \int \, r^k \, \vert f (\bm{r}) \vert^2 \,
  \mathrm{d}^3 \bm{r} < \infty \Bigr\} \, ,
  \qquad k = -1, 1, 2, \dots \; .
\end{equation}
For $k=0$, we retrieve the Hilbert space $L^{2} (\mathbb{R}^{3})$ of
square integrable functions defined by (\ref{HilbertL^2}).

As long as the principal quantum number $N$ is not too negative, a
Slater-type function $\chi_{N, L}^{M} (\beta, \bm{r})$ is for a fixed
value of $k$ an element of the weighted Hilbert space $L_{r^k}^{2}
(\mathbb{R}^{3})$. In this case, Guseinov's approach
\cite{Guseinov/1980a,Guseinov/2001a,Guseinov/2002c}, who constructed
one-range addition theorems by expanding $\chi_{N, L}^{M} (\beta, \bm{r}
\pm \bm{r}')$ in terms of his complete and orthonormal functions $\{
\prescript{}{k}{\Psi}_{n, \ell}^{m} (\gamma, \bm{r}) \}_{n, \ell, m}$
with in general different scaling parameters $\beta \neq \gamma > 0$, is
mathematically sound. For fixed $k = -1, 0, 1, 2, \dots$, Guseinov
constructed expansions that converge in the mean with respect to the norm of
$L_{r^k}^{2} (\mathbb{R}^3)$.

However, Guseinov replaced in his one-range addition theorems his
complete and orthonormal functions $\prescript{}{k}{\Psi}_{n, \ell}^{m}
(\beta, \bm{r})$ by nonorthogonal Slater-type functions with integral
principal quantum numbers via \cite[Eq.\ (6.4)]{Weniger/2007b}
\begin{align}
  \label{GusFun2STF}
  \prescript{}{k}{\Psi}_{n, \ell}^{m} (\beta,
  \bm{r}) & \; = \; 2^{\ell} \, \left[ \frac{(2\beta)^{k+3} \,
      (n+\ell+k+1)!}{(n-\ell-1)!} \right]^{1/2}
  \notag \\
  & \qquad \times \, \sum_{\nu=0}^{n-\ell-1} \,
  \frac{(-n+\ell+1)_{\nu} \, 2^{\nu}}{(2\ell+k+\nu+2)! \, \nu!} \,
  \chi_{\nu+\ell+1, \ell}^{m} (\beta, \bm{r}) \, .
\end{align}
This is still legitimate. However, Guseinov also rearranged the order of
summations of the resulting expansions. In this way, Guseinov formally
constructed expansions of Slater-type functions $\chi_{N, L}^{M} (\beta,
\mathbf{r} \pm \mathbf{r}')$ with in general nonintegral principal
quantum numbers $N \in \mathbb{R} \setminus \mathbb{N}$ in terms of
Slater-type functions $\chi_{n, \ell}^{m} (\beta, \mathbf{r})$ with
integral principal quantum numbers $n \in \mathbb{N}$ located at a
different center (see also \cite[Section 6]{Weniger/2007b}).

Slater-type functions are complete in all Hilbert spaces, which Guseinov
implicitly used (for an explicit proof, see \cite[Section
4]{Klahn/Bingel/1977b}), but not orthogonal. Thus, Guseinov's approach
corresponds to the transformation of an expansion in terms of a complete
and orthogonal function set to an expansion in terms a of a complete, but
nonorthogonal function set.

Unfortunately, the completeness of a nonorthogonal function set in a
Hilbert space does not suffice to guarantee that an essentially arbitrary
element of this Hilbert space can be expanded in terms of this function
set (the nonanalyticity of certain Laguerre series discussed in this
article is just another confirmation of a much more general fact). This
insufficiency is well documented both in the mathematical literature (see
for example \cite[Theorem 10 on p.\ 54]{Davis/1989} or \cite[Section
1.4]{Higgins/1977}) as well as in the literature on electronic structure
calculations
\cite{Klahn/1975,Klahn/1981,Klahn/Bingel/1977a,Klahn/Bingel/1977b,%
  Klahn/Bingel/1977c,Klahn/Morgan/1984}), but nevertheless frequently
overlooked.  Horrifying examples of nonorthogonal expansions with
pathological properties can be found in \cite[Section III.I]{Klahn/1981}.

Consequently, it is not at all obvious whether the mathematical
manipulations, that produced Guseinov's rearranged one-range addition
theorems, are legitimate and lead to expansions that are mathematically
meaningful. This has to be checked. So far, Guseinov has categorically
denied that there might be any problem with the legitimacy of his
rearrangements \cite[pp.\ 8 and 23 - 24]{Guseinov/2007a}.

Since Guseinov's original addition theorems are expansions in terms of
generalized Laguerre polynomials, and since the radial parts of
Slater-type functions are after the cancellation of common exponentials
nothing but powers, the results of this article about the analyticity of
Laguerre series could in principle be used to investigate whether
Guseinov's rearrangements are legitimate. For that purpose, it would be
necessary to determine the decay rates and the sign patterns of the
coefficients of the generalized Laguerre polynomials occurring in
Guseinov's original addition theorems.

Unfortunately, one-range addition theorems for Slater-type functions are
fairly complicated mathematical objects, and the coefficients of the
generalized Laguerre polynomials are according to
(\ref{OneRangeAddTheor_w}) essentially three-dimensional overlap
integrals. Thus, we would be confronted with enormous and possibly even
unsurmountable technical problems if we try to analyze the decay rates
and sign patters of the coefficients of these Laguerre series.

Fortunately, some insight can be gained by analyzing not the
comparatively complicated one-range addition theorems, but their much
simpler one-center limits. These are also expansions in terms of
generalized Laguerre polynomials, albeit with much simpler coefficients.
It will become clear later that this approach cannot answer all questions
in interest. Nevertheless, it is better than nothing since we can obtain
at least some nontrivial answers.

Let us now consider the one-center expansion of a Slater-type function
$\chi_{N, L}^{M} (\beta, \mathbf{r})$ with in general nonintegral
principal quantum numbers $N \in \mathbb{R} \setminus \mathbb{N}$ in
terms of Guseinov's functions $\prescript{}{k}{\Psi}_{n, \ell}^{m}
(\beta, \bm{r})$ with \emph{equal} scaling parameters $\beta > 0$
\cite[Eq.\ (5.7)]{Weniger/2007c}:
\begin{align}
  \label{NISTF2Gusfun_EqScaPar}
  & \chi_{N, L}^{M} (\beta, \bm{r}) \; = \; \frac
  {(2\gamma)^{-(k+3)/2}}{2^{N-1}} \, \Gamma (N+L+k+2)
  \notag \\
  & \qquad \times \, \sum_{\nu=0}^{\infty} \,
  \frac{(-N+L+1)_{\nu}}{\bigl[ (\nu+2L+k+2)! \, \nu! \bigr]^{1/2}}
  \, \prescript{}{k}{\Psi}_{\nu+L+1, L}^{M} (\beta, \bm{r}) \, ,
  \notag \\
  & \qquad \qquad N \in \mathbb{R} \setminus \mathbb{N} \, ,
    \qquad \beta > 0 \, , \qquad k = -1, 0, 1, 2, \dots \, .
\end{align}
If $N \in \mathbb{N}$ and $N \ge L+1$, the infinite series on the
right-hand side terminates because of the Pochhammer symbol
$(-N+L+1)_{\nu}$.

The expansion (\ref{NISTF2Gusfun_EqScaPar}), which corresponds to the
one-center limit $\bm{r}' = \bm{0}$ of Guseinov's one-range addition
theorem for $\chi_{N, L}^{M} (\beta, \mathbf{r} \pm \bm{r}')$ with equal
scaling parameters, is a special case of the Laguerre series
(\ref{GenPow2GLag}) for $z^{\rho}$ with in general nonintegral $\rho \in
\mathbb{R} \setminus \mathbb{N}_0$. As discussed in Section
\ref{Sec:LagSerGeneralPowerFunction}, a rearrangement of the Laguerre
series (\ref{GenPow2GLag}) for $z^{\rho}$ is legitimate if and only if
$\rho$ is a nonnegative integer, i.e., if $\rho = m \in \mathbb{N}_{0}$.
If $\rho \notin \mathbb{N}_{0}$, $z^{\rho}$ is not analytic at $z=0$ and
we formally obtain the power series (\ref{Chk_x^m_GlagPol_6}), which is
not a mathematically meaningful object since an infinite number of its
power series coefficients are according to
(\ref{Lim_Ser_Chk_x^m_GlagPol_6}) infinite in magnitude.

Thus, in the case of equal scaling parameters $\beta > 0$, the one-center
limit $\bm{r}' = \bm{0}$ of Guseinov's rearranged addition theorem for
Slater-type functions $\chi_{N, L}^{M} (\beta, \bm{r} \pm \bm{r}')$ does
not exist if the principal quantum number $N$ is nonintegral, i.e., if $N
\in \mathbb{R} \setminus \mathbb{N}$.

Let us now consider the one-center expansion of a Slater-type function
$\chi_{N, L}^{M} (\beta, \mathbf{r})$ with in general nonintegral
principal quantum numbers $N \in \mathbb{R} \setminus \mathbb{N}$ in
terms of Guseinov's functions $\prescript{}{k}{\Psi}_{n, \ell}^{m}
(\gamma, \bm{r})$ with \emph{different} scaling parameters $\beta \ne
\gamma > 0$ \cite[Eq.\ (5.9)]{Weniger/2007c}:
\begin{align}
  \label{Expand_NISTF2Gusfun_DiffScaPar_1}
  \chi_{N, L}^{M} (\beta, \bm{r}) & \; = \; \frac
  {(2\gamma)^{L+(k+3)/2} \, \beta^{N-1}}{[\beta+\gamma]^{N+L+k+2}}
  \, \frac{\Gamma (N+L+k+2)}{(2L+k+2)!}
  \notag \\
  & \qquad \times \sum_{\nu=0}^{\infty} \, \left[
    \frac{(\nu+2L+k+2)!}{\nu!} \right]^{1/2} \,
  \prescript{}{k}{\Psi}_{\nu+L+1, L}^{M} (\gamma, \bm{r})
  \notag \\
  & \qquad \qquad \times {}_2 F_1 \left(-\nu, N+L+k+2; 2L+k+3;
    \frac{2\gamma}{\beta+\gamma} \right) \, .
\end{align}
The expansion (\ref{Expand_NISTF2Gusfun_DiffScaPar_1}), which corresponds
to the one-center limit $\bm{r}' = \bm{0}$ of Guseinov's one-range
addition theorem for $\chi_{N, L}^{M} (\beta, \mathbf{r} \pm \bm{r}')$
with different scaling parameters, is a special case of the following
Laguerre series \cite[Eq.\ (6.12)]{Weniger/2007b}:
\begin{align}
  \label{ExpoPow2GLag}
  z^{\rho} \, \mathrm{e}^{u z} & \; = \; (1-u)^{-\alpha-\rho-1} \,
  \frac{\Gamma (\alpha+\rho+1)}{\Gamma (\alpha+1)}
  \notag \\
  & \qquad \times \, \sum_{n=0}^{\infty} \, {}_2 F_1 \left(-n,
    \alpha+\rho+1; \alpha+1; \frac{1}{1-u} \right) \, L_{n}^{(\alpha)}
  (z) \, ,
  \notag \\
  & \qquad \qquad \rho \in \mathbb{R} \setminus \mathbb{N}_{0} \, ,
  \quad \Re (\rho+\alpha) > - 1 \, , \! \quad u \in (-\infty, 1/2) \, .
\end{align}
This expansion can be derived with the help of (\ref{Gr_7.414.7}). The
condition $-\infty < u < 1/2$ is necessary to guarantee its convergence
in the mean with respect to the norm (\ref{Def:Norm_Lag}) of the weighted
Hilbert space $L^{2}_{z^{\alpha} \mathrm{e}^{-z}} \bigl( [0,
\infty)\bigr)$. For $u=0$, the terminating Gaussian hypergeometric series
${}_{2} F_{1}$ can be expressed in closed form with the help of Gauss'
summation theorem \cite[p.\ 40]{Magnus/Oberhettinger/Soni/1966}, and
(\ref{ExpoPow2GLag}) simplifies to give (\ref{GenPow2GLag}).

If we insert the explicit expression (\ref{GLag_1F1}) of the generalized
Laguerre polynomials into (\ref{ExpoPow2GLag}) and interchange the order
of summations, we also obtain a formal power series in $z$.
Unfortunately, an analysis of the resulting power series becomes very
difficult because of the terminating Gaussian hypergeometric series ${}_2
F_1$ in (\ref{ExpoPow2GLag}). An analysis of the behavior of this ${}_2
F_1$ as $n \to \infty$ would most likely be a nontrivial research project
in its own right. However, we can argue that $z^{\rho} \exp (u z)$ is
only analytic at $z=0$ if $\rho = m \in \mathbb{N}_0$, yielding the
expansion $z^m \exp (u z) = \sum_{n=0}^{\infty} u^n z^{m+n}/n!$. If
$\rho$ is nonintegral, a power series expansion of $z^{\rho} \exp (u z)$
about $z=0$ does not exist.

Thus, also for different scaling parameters $\beta \neq \gamma$, the
one-center limit $\bm{r}' = \bm{0}$ of the rearranged addition theorems
for $\chi_{N, L}^{M} (\beta, \bm{r} \pm \bm{r}')$ does not exist if the
principal quantum number $N$ is nonintegral, i.e., if $N \in \mathbb{R}
\setminus \mathbb{N}$.

These observations are quite consequential: Rearranged one-range addition
theorems for Slater-type functions $\chi_{N, L}^{M} (\beta, \bm{r} \pm
\bm{r}')$ with nonintegral principal quantum numbers $N \neq 1, 2, \dots$
play a central role in numerous articles by Guseinov and coworkers on
one-range addition theorems
\cite{Guseinov/2002b,Guseinov/2002d,Guseinov/2003b,Guseinov/2003d,%
  Guseinov/2003e,Guseinov/2004a,Guseinov/2004b,Guseinov/2004c,%
  Guseinov/2004d,Guseinov/2004i,Guseinov/2004k,%
  Guseinov/2005a,Guseinov/2005c,Guseinov/2005g,Guseinov/2006a,%
  Guseinov/Mamedov/2002d,Guseinov/Mamedov/2004b,Guseinov/Mamedov/2004d,%
  Guseinov/Mamedov/2004e,Guseinov/Mamedov/2005c,%
  Guseinov/Mamedov/Suenel/2002}.

By analyzing the rearrangement of the Laguerre series of the
comparatively simple functions $z^{\rho}$ and $z^{\rho} \exp (u z)$, we
could arrive at some conclusions about the legitimacy of Guseinov's
approach. Nevertheless, some interesting questions are still open. For
example, the nonanalyticity arguments presented here allow no conclusions
about the validity of Guseinov's rearrangements in the case of
Slater-type functions with integral principal quantum numbers.

Another interesting but open question is whether Guseinov's
rearrangements produce in the case of nonintegral principal quantum
numbers one-range addition theorems that are invalid for the whole
argument set $\mathbb{R}^{3} \times \mathbb{R}^{3}$, or whether only the
one-center limits of these addition theorems are invalid. This is a
practically very relevant question. If only the one-center limits are
invalid, then it would be conceivable that Guseinov's rearranged
one-range addition theorems might be mathematically meaningful or
possibly even numerically useful in a restricted sense as approximations,
although they do not exist for the whole argument set $\mathbb{R}^3
\times \mathbb{R}^3$. This remains to be investigated. But the burden of
proof lies in all cases with Guseinov.

\typeout{==> Section: Summary and Conclusions}
\section{Summary and Conclusions}
\label{Sec:SummaryConclusions}

The generalized Laguerre polynomials belong to the so-called classical
orthogonal polynomials of mathematical physics, and they are
characterized by the orthogonality relationship
(\ref{GLag_Orthogonality}) that involves an integration over the
semi-infinite positive real axis. Accordingly, generalized Laguerre
polynomials can be used for the representation of functions on unbounded
domains, and in particular also for the representation of the radial
parts of functions $f \colon \mathbb{R}^{3} \to \mathbb{C}$ expressed in
terms of spherical polar coordinates.

It is generally accepted that expansions in terms of orthogonal
polynomials have many highly advantageous features. However, there is one
undeniable drawback: Normally, expansions in terms of orthogonal
polynomials converge in the mean with respect to the norm of the
corresponding Hilbert space, but not necessarily pointwise. Accordingly,
orthogonal expansions are not necessarily the best choice if the
\emph{local} properties of a function matter. As documented by the
popularity of Pad\'{e} approximants, power series are also not free of
weaknesses, but at least in the vicinity of the expansion point, power
series are normally very convenient and very useful for an accurate
description of the local properties of a function. Therefore, the
construction of power series expansions for functions defined by
orthogonal expansions should be of principal interest.

This article describes a complementary treatment of Laguerre series of
the type of (\ref{f_Exp_GLag}). Normally, one starts from a \emph{known}
function $f (z)$ belonging to the weighted Hilbert space
$L^{2}_{z^{\alpha} \mathrm{e}^{-z}} \bigl( [0, \infty)\bigr)$ defined by
(\ref{HilbertL^2_Lag}), and one tries to determine the coefficients
$\lambda_{n}^{(\alpha)}$ via (\ref{f_Exp_GLag_b}) by exploiting the
orthogonality of the generalized Laguerre polynomials.

In this article, it is instead assumed that only the Laguerre series
coefficients $\lambda_{n}^{(\alpha)}$ are known, either in the form of
explicit expressions or numerically, but not the function $f (z)$. With
the help of the transformation formula (\ref{Rearr_f_Exp_GLag}), it is
then possible to construct a \emph{formal} power series expansion of the
\emph{unknown} function represented by the Laguerre series.

This approach does not guarantee success since there are many functions
which belong to the Hilbert space $L^{2}_{z^{\alpha} \mathrm{e}^{-z}}
\bigl( [0, \infty)\bigr)$ but which are not analytic at the origin. A
simple example of such a nonanalytic function possessing a Laguerre
series is the power function $z^{\rho}$ with nonintegral $\rho \in
\mathbb{R} \setminus \mathbb{N}_{0}$. As discussed in Section
\ref{Sec:LagSerGeneralPowerFunction}, it is possible to construct a
formal power series for $z^{\rho}$ from its Laguerre series
(\ref{GenPow2GLag}), but the resulting power series is mathematically
meaningless since it contains infinitely many power series coefficients
that are infinite in magnitude.

Thus, the key question is whether the power series for the unknown
function obtained via (\ref{Rearr_f_Exp_GLag}) is mathematically
meaningless, or whether this power series represents an analytic function
in the sense of complex analysis. It seems that this question has not
been treated properly in the literature yet. I am only aware of short
remarks by Gottlieb and Orszag \cite[p.\ 42]{Gottlieb/Orszag/1977} and by
Doha \cite[p.\ 5452]{Doha/2003}, respectively, who had stated that a
Laguerre series of the type of (\ref{f_Exp_GLag}) converges faster than
algebraically if the function under consideration is analytic at the
origin. But this statement is imprecise and ignores the pivotal role
played by divergent, but summable inner $\mu$ series in
(\ref{Rearr_f_Exp_GLag}). Some general aspects of the summation of
divergent series are reviewed in Appendix \ref{App:DivergentSeries}.

By analyzing the convergence properties of the inner $\mu$ series in
(\ref{Rearr_f_Exp_GLag}), some simple \emph{sufficient} conditions can be
formulated which guarantee that the resulting power series is
mathematically meaningful and represents an analytic function.

As discussed in Section
\ref{Sec:ExponentiallyFactoriallyDecayingSeriesCoefficients}, the most
benign situation occurs if the Laguerre series coefficients
$\lambda_{n}^{(\alpha)}$ decay exponentially or even factorially as $n
\to \infty$. Then, the inner $\mu$ series in (\ref{Rearr_f_Exp_GLag})
converge and the resulting power series is mathematically meaningful and
represents an analytic function. In this way, numerous generating
functions for the generalized Laguerre polynomials can be rederived
easily.

As discussed in Section
\ref{Sec:AlgebraicallyDecayingSeriesCoefficients}, a much more
interesting situation occurs if the coefficients $\lambda_{n}^{(\alpha)}$
decay algebraically in magnitude as $n \to \infty$. If algebraically
decaying coefficients $\lambda_{n}^{(\alpha)}$ ultimately have the same
sign, the inner $\mu$ series in (\ref{Rearr_f_Exp_GLag}) diverge. This
alone wound not necessarily be such a bad thing, but it is not possible
to sum these divergent series to something finite. Accordingly, the
transformation formula (\ref{Rearr_f_Exp_GLag}) leads to a power series
expansion having infinitely many series coefficients that are infinite in
magnitude. This simply means that a power series of the type of
(\ref{PowSer_f}) does not exist because the function under consideration
is not analytic at the origin.  An example is the Laguerre series
(\ref{GenPow2GLag}) for $z^{\rho}$: Its series coefficients decay
algebraically in magnitude and ultimately all have the same sign. This
observation suffices to show once more that $z^{\rho}$ is not analytic at
the origin if $\rho \in \mathbb{R} \setminus \mathbb{N}_{0}$.

The situation changes radically if the coefficients
$\lambda_{n}^{(\alpha)}$ decay algebraically in magnitude as $n \to
\infty$, but ultimately have \emph{strictly alternating} signs. Then, the
inner $\mu$ series in (\ref{Rearr_f_Exp_GLag}) still do not converge, but
now summation techniques can be used to associate finite values to these
divergent alternating series. In such a case, the formal power series
obtained via (\ref{Rearr_f_Exp_GLag}) is mathematically meaningful and
represents an analytic function.

For example, the additional sign factor $(-1)^{n}$ introduced into the
Laguerre series (\ref{GenPow2GLag}) for $z^{\rho}$ yields the Laguerre
series (\ref{LagSer_AltSerCoeffs}). The application of
(\ref{Rearr_f_Exp_GLag}) to (\ref{LagSer_AltSerCoeffs}) leads to
divergent inner $\mu$ series which correspond to divergent hypergeometric
series ${}_{1} F_{0}$. But these series can be summed by analytic
continuation since they are special cases of the binomial series
(\ref{BinomSer}). Thus, the transformation formula
(\ref{Rearr_f_Exp_GLag}) ultimately produces according to
(\ref{Rearr_G_Exp_GLag_4}) a confluent hypergeometric series ${}_{1}
F_{1}$, which is an analytic function in every neighborhood of the origin
and which is also a known generating function of the generalized Laguerre
polynomials.

The summability approach pursued in Section
\ref{Sec:AlgebraicallyDecayingSeriesCoefficients} can also be used in the
case of Laguerre series with more complicated coefficients. An example is
the Laguerre series (\ref{LagSer_2H1_a_b_c_alpha}). The application of
the transformation formula (\ref{Rearr_f_Exp_GLag}) leads to divergent,
but summable $\mu$ series that can be expressed as a divergent Gaussian
hypergeometric series ${}_{2} F_{1}$. But again, the summation of this
divergent series is almost trivially simple because many convenient
analytic continuation formulas for a ${}_{2} F_{1}$ are known.

Thus, the summation of a divergent inner $\mu$ series can be accomplished
by explicit analytic continuation formulas if it can be expressed either
as a divergent binomial series ${}_{1} F_{0}$ or as a divergent Gaussian
hypergeometric series ${}_{2} F_{1}$. Unfortunately, this is no longer
possible if the $\mu$ series corresponds to a divergent generalized
hypergeometric series ${}_{p+1} F_{p}$ with $p \ge 2$. As discussed in
Section \ref{Sec:ComputationalApproaches}, analytic continuation formulas
for more complicated generalized hypergeometric series are either not
known at all or at least much more complicated than the corresponding
formulas for Gaussian hypergeometric series.

Thus, we can only hope to find convenient analytical expressions for
power series coefficients $\gamma_{n}$ if the corresponding algebraically
decaying and ultimately strictly alternating Laguerre series coefficients
$\lambda_{n}^{(\alpha)}$ possess a very simple structure.  But this is a
typical limitation of all analytical manipulations. As a viable
alternative, we can try to use instead techniques that accomplish a
summation of divergent inner $\mu$ series by purely \emph{numerical}
means. Such an approach has the additional advantage that it can be
applied if only the numerical values of a finite set of Laguerre series
coefficients $\lambda_{n}^{(\alpha)}$ are available.

In Section \ref{Sec:ComputationalApproaches}, it is shown that certain
nonlinear sequence transformations, whose properties are reviewed in
Appendix \ref{App:SequenceTransformations}, are indeed able to sums
divergent alternating hypergeometric series ${}_{p+1} F_{p}$ with $p \ge
1$ quite effectively. It may be surprising for nonspecialists that the
best summation results were not obtained by Wynn's celebrated epsilon
algorithm (\ref{eps_al}), which produces Pad\'{e} approximants if the
input data are the partial sums of a power series, but by the Levin-type
transformation (\ref{dWenTr}).

In Section \ref{Sec:GuseinovsRearrangedOne-RangeAdditionTheorems}, the
legitimacy of the rearrangement of Guseinov's one-range addition theorems
for Slater-type functions with in general nonintegral principal quantum
numbers is analyzed. Originally, Guseinov
\cite{Guseinov/1980a,Guseinov/2001a,Guseinov/2002c} had derived addition
theorems for Slater-type functions by expanding them in terms of the
complete and orthonormal functions (\ref{Def_Psi_Guseinov}) whose radial
parts consist of generalized Laguerre polynomials. In the next step,
Guseinov replaced the functions (\ref{Def_Psi_Guseinov}) according to
(\ref{GusFun2STF}) by Slater-type functions with integral principal
quantum numbers, which are complete, but not orthogonal, and he also
rearranged the order of the resulting summations. In this way, Guseinov
essentially replaced expansions in terms of generalized Laguerre
polynomials by power series expansions. Thus, the results of this article
can in principle be used to check the legitimacy of Guseinov's
manipulations.

Since, however, one-range addition theorems for Slater-type functions are
fairly complicated mathematical objects, the necessary determination of
the decay rates and the sign patterns of the coefficients of the
generalized Laguerre polynomials is very difficult. Fortunately, one can
gain at least some insight by analyzing not the complicated addition
theorems, but their much simpler one-center limits. In this way, it can
for instance be shown that Guseinov's rearranged addition theorems for
Slater-type functions with nonintegral principal quantum numbers do not
exist for the whole argument range $\mathbb{R}^{3} \times
\mathbb{R}^{3}$.


\begin{appendix}
\typeout{==> Appendix A: Divergent Series}
\section{Divergent Series}
\label{App:DivergentSeries}

Divergent series have been a highly controversial topic that played a
major role in the development of mathematical analysis (see for example
the article by Burkhardt \cite{Burkhardt/1911} or the very recent book by
Ferraro \cite{Ferraro/2008}), and even now there is still a lot of active
research on divergent series going on (see for example the recent
monographs by Balser \cite{Balser/1994,Balser/1999}, Boutet de Monvel
\cite{BoutetDeMonvel/1994}, Candelpergher, Nosmas, and Pham
\cite{Candelpergher/Nosmas/Pham/1993}, and Sternin and Shatalov
\cite{Sternin/Shatalov/1996}, or the review by Delabaere and Pham
\cite{Delabaere/Pham/1999}).

As for instance discussed in articles by Barbeau \cite{Barbeau/1979},
Barbeau and Leah \cite{Barbeau/Leah/1976}, Kozlov \cite{Kozlov/2007}, and
Varadarajan \cite{Varadarajan/2007}, already Euler had frequently used
divergent series. Later, when the concept of convergence was better
understood, Euler's admittedly somewhat informal treatment of divergent
series was criticized, and a strong tendency emerged to ban divergent
series completely from the realm of rigorous mathematics. This criticism
culminated in Abel's famous quotation from the year 1828, which very well
expressed the prevalent attitude of most mathematicians during a large
part of the nineteenth century and which can for instance be found in
Littlewood's preface of Hardy's posthumously published classic on
divergent series \cite{Hardy/1949}:
\begin{quote}
  \sl Divergent series are the invention of the devil, and it is shameful
  to base on them any demonstration whatsoever.
\end{quote}
Although this disdain of divergent series was a temporary phenomenon of
the nineteenth century, its consequences are nevertheless felt today. As
emphasized by Suslov \cite[p.\ 1191]{Suslov/2005}, the standard
university curricula in mathematical analysis were formulated in the
middle of the nineteenth century. But this was a time when divergent
series were wrongly considered to be essentially an aberration of the
pre-rigorous mathematical past. As a consequence, divergent series and
their summation are not part of the regular training of mathematicians
and theoretical physicists, which in my opinion is totally inappropriate.
In this context, a short article by Rubel \cite{Rubel/1989} may also be
of interest.

At the end of the nineteenth century it was clear that the attempts of
mathematical orthodoxy to reject divergent series as unfounded had
failed. Firstly, divergent series turned out to be too useful to be
abandoned. For example, many special functions possess so-called
asymptotic expansions which normally are factorially divergent inverse
power series. In spite of their divergence, suitably truncated asymptotic
series can provide excellent approximations at least for sufficiently
large arguments.

Secondly, the work of mathematicians like
\begin{itemize}
\item Poincar\'{e}, whose work in astronomy inspired him to formulate a
  mathematically rigorous theory of asymptotic series that typically
  diverge factorially,
\item Borel, who showed that factorially divergent series can be summed
  via Laplace-type integral representations,
\item Pad\'{e}, who introduced his celebrated rational approximants that
  are often able to sum divergent power series,
\item Stieltjes, who showed that certain divergent series can be
  identified with continued fractions,
\end{itemize}
ultimately led in the later part of the nineteenth century to a theory
which made it possible to use divergent series in a mathematically
rigorous way. Even more important from a practical point of view is that
their work showed that divergent series can actually be used for
computational purposes if they are combined with suitable summation
techniques.

In physics, divergent series are now indispensable. It is generally
accepted that perturbation theory is the most important \emph{systematic}
approximation procedure in theoretical physics. But already in 1952,
Dyson \cite{Dyson/1952} had argued that perturbation expansions in
quantum electrodynamics must diverge factorially, and since the seminal
work of Bender and Wu on the perturbation expansions of anharmonic
oscillators \cite{Bender/Wu/1969,Bender/Wu/1971} it has been clear that
quantum mechanical perturbation theory produces almost by default
factorially divergent perturbation expansions. A good source on divergent
perturbation expansions in quantum mechanics and in higher field theories
is the book edited by Le Guillou and Zinn-Justin
\cite{LeGuillou/Zinn-Justin/1990} where many of the relevant articles are
reprinted.

Even if we agree that divergent series are now indispensable in physics,
we nevertheless cannot expect that rigorously minded mathematicians are
necessarily satisfied with the way divergent series are typically used in
physics. Quite instructive is the following remark by Haldane from the
year 1941, which is quoted in a book by K\"{o}rner \cite[p.\
426]{Koerner/1988}:
\begin{quote}
  \sl Cambridge is full of mathematicians who have been so corrupted by
  quantum mechanics that they use series which are clearly divergent, and
  not even proved to be summable.
\end{quote}
I think that the criticism, which is implicit in Haldane's remark, cannot
be dismissed lightly. Let us for instance assume that we want to express
a physical quantity by a divergent series, which for example may be a
factorially divergent Rayleigh-Schr\"{o}dinger perturbation expansion.
Since divergent series are \emph{a priori} mathematically undefined, we
have to show that it is indeed possible to associate something finite --
the value of the physical quantity -- to the divergent series in a
mathematically meaningful way, or to put it differently, we have to show
that the divergent series is summable to the correct finite result by an
appropriate summation method.

Unfortunately, it is often extremely difficult to prove\ this rigorously.
In such a situation, physicists tend to rely on their intuition and do
not bother to try to formulate difficult proofs. Quite often, such a
pragmatic approach is remarkably successful, but it should also be clear
that mathematicians do not necessarily like that. Intuition can be
misleading.  Therefore, occasional unpleasant surprises and even
catastrophic failures cannot be ruled out.

\typeout{==> Appendix B: Sequence Transformations}
\section{Sequence Transformations}
\label{App:SequenceTransformations}

The fact that appropriate summation methods make it possible to use
divergent series for computational purposes raises the question which of
the numerous known summation techniques are best suited for the numerical
evaluation of the divergent series occurring in this article.  Typically,
we are confronted with alternating series whose terms grow in magnitude
like a fixed power of the index. Based on my own practical experience, I
propose to use so-called \emph{sequence transformations}, which are
purely numerical techniques to transform a slowly convergent or divergent
sequence $\{ s_n \}_{n=0}^{\infty}$ to another sequence $\{
s_{n}^{\prime} \}_{n=0}^{\infty}$ with hopefully better convergence
properties.

For those interested in the history of sequence transformations, I
recommend a monograph by Brezinski \cite{Brezinski/1991a}, which
discusses earlier work starting from the 17th century until 1945, as well
as two articles by Brezinski \cite{Brezinski/1996,Brezinski/2000c}, which
emphasize more recent developments.

The basic assumption of all sequence transformations is that the elements
of a slowly convergent or divergent sequence $\{ s_n \}_{n=0}^{\infty}$,
which could be the partial sums $s_n = \sum_{k=0}^{n} a_{k}$ of an
infinite series, can for all indices $n$ be partitioned into a
(generalized) limit $s$ and a remainder or truncation error $r_n$
according to
\begin{equation}
  \label{s_n_r_n}
  s_n \; = \; s + r_n \, , \qquad n \in \mathbb{N}_0 \, .
\end{equation}
If the sequence $\{ s_n \}_{n=0}^{\infty}$ converges to $s$, the
remainders $r_n$ in (\ref{s_n_r_n}) can be made negligible by increasing
$n$ as much as necessary. But many sequences converge so slowly that this
is not feasible. Increasing the index $n$ also does not help in the case
of the divergent series.

Alternatively, one can try to improve convergence or accomplish a
summation by computing approximations to the remainders $r_{n}$ which are
then eliminated from the sequence elements $s_{n}$. At least
conceptually, this is what a sequence transformation tries to do.

With the exception of a few practically more or less irrelevant model
problems, sequence transformations can only eliminate approximations to
the remainders. Thus, the elements of the transformed sequence $\{
s_n^{\prime} \}_{n=0}^{\infty}$ are also be of the type of
(\ref{s_n_r_n}), which means that a transformed sequence element
$s^{\prime}_n$ can also be partitioned into the same (generalized) limit
$s$ and a transformed remainder $r_n^{\prime}$ according to
\begin{equation}
  s_n^{\prime} \; = \; s + r_n^{\prime} \, ,
  \qquad n \in \mathbb{N}_0 \, .
\end{equation}
The transformed remainders $\{ r_n^{\prime} \}_{n=0}^{\infty}$ are
normally different from zero for all finite values of $n$. However,
convergence is accelerated if the transformed remainders $\{ r_n^{\prime}
\}_{n=0}^{\infty}$ vanish more rapidly than the original remainders $\{
r_n \}_{n=0}^{\infty}$, and a divergent sequence is summed if the
transformed remainders vanish at all as $n \to \infty$.

Before the invention of electronic computers, mainly \emph{linear}
sequence transformations were used, which compute the elements of the
transformed sequence $\{ s_n^{\prime} \}_{n=0}^{\infty}$ as weighted
averages of the elements of the input sequence $\{ s_n
\}_{n=0}^{\infty}$ according to
\begin{equation}
s_n^{\prime} \; = \; \sum_{k=0}^{n} \, \mu_{n k} \, s_k \, .
\label{Mat_tr}
\end{equation}
The theoretical properties of these matrix transformations are now very
well understood and discussed in books by Hardy \cite{Hardy/1949}, Knopp
\cite{Knopp/1964}, Petersen \cite{Petersen/1966}, Peyerimhoff
\cite{Peyerimhoff/1969}, Zeller and Beekmann \cite{Zeller/Beekmann/1970},
Powell and Shah, and Boos \cite{Boos/2000}. Their main appeal lies in the
fact that based on the work of Toeplitz \cite{Toeplitz/1911} some
necessary and sufficient conditions for the weights ${\mu}_{n k}$ in
(\ref{Mat_tr}) could be formulated which guarantee that the application
of such a matrix transformation to a convergent sequence $\{ s_n
\}_{n=0}^{\infty}$ yields a transformed sequence $\{ s_n^{\prime}
\}_{n=0}^{\infty}$ converging to the same limit $s = s_{\infty}$.

There are also some so-called Tauberian theorems which show rigorously
that certain orders of for example H\"{o}lder's or Ces\`{a}ro's summation
method are needed to sum a divergent series whose terms grow in magnitude
like a fixed power of the index \cite[Theorem 39 on p.\ 95 and Theorem 49
on p.\ 103]{Hardy/1949}. Since the terms of the divergent inner $\mu$
series considered in Sections
\ref{Sec:AlgebraicallyDecayingSeriesCoefficients} and
\ref{Sec:ComputationalApproaches} all grow in magnitude like a fixed
power of the index, it looks like a natural idea to use either
H\"{o}lder's or Ces\`{a}ro's summation method for the summation of these
divergent series. Unfortunately, the situation is not so simple. Since we
also want to accomplish something useful if we only know a comparatively
small number of numerically determined Laguerre series coefficients
$\lambda_{n}^{(\alpha)}$, we should focus our attention on those purely
numerical summation techniques that promise to be particularly efficient.

From a purely theoretical point of view, regularity is extremely
desirable and greatly facilitates the formulation of nice mathematical
proofs, but from a practical point of view, regularity is a serious
disadvantage. This probably sounds paradoxical. However, Wimp remarks in
the preface of his book \cite[p.\ X]{Wimp/1981} that the size of the
domain of regularity of a transformation and its efficiency seem to be
inversely related.  Accordingly, regular matrix transformations are in
general at most moderately powerful, and the popularity of most linear
transformations as computational tools has declined considerably in
recent years.

Nonlinear sequence transformations have largely complementary properties:
They are nonregular, which means that the convergence of the transformed
sequence is not guaranteed, let alone to the correct limit. In addition,
their theoretical properties are far from being completely understood.
Thus, from a purely theoretical point of view, nonlinear sequence
transformations have many disadvantages, but they often accomplish
spectacular transformation results which are clearly beyond the reach of
regular matrix transformations. Consequently, nonlinear transformations
now clearly dominate practical applications and -- albeit to a lesser
extend -- also theoretical work. Detailed treatments of their theoretical
properties and long lists of successful applications can be found in
monographs by Brezinski
\cite{Brezinski/1977,Brezinski/1978,Brezinski/1980a}, Brezinski and
Redivo Zaglia \cite{Brezinski/RedivoZaglia/1991a}, Cuyt
\cite{Cuyt/1988,Cuyt/1994}, Cuyt and Wuytack \cite{Cuyt/Wuytack/1987},
Delahaye \cite{Delahaye/1988}, Liem, L\"{u}, and Shih
\cite{Liem/Lu/Shih/1995}, Marchuk and Shaidurov
\cite{Marchuk/Shaidurov/1983}, Sidi \cite{Sidi/2003}, Walz
\cite{Walz/1996}, and Wimp \cite{Wimp/1981}, or in reviews by Caliceti,
Meyer-Hermann, Ribeca, Surzhykov, and Jentschura
\cite{Caliceti/Meyer-Hermann/Ribeca/Surzhykov/Jentschura/2007}, Homeier
\cite{Homeier/2000a}, and myself \cite{Weniger/1989}).

In spite of their undeniable usefulness, nonlinear sequence
transformations do not necessarily get the attention they deserve, in
particular in the more theoretically oriented mathematical literature.
For example, a (very) condensed review of the classical linear
summability methods associated with the names of Ces\`{a}ro, Abel, and
Riesz can be found in Zayed's relatively recent book \cite[Chapter
1.11.1]{Zayed/1996}, but the more powerful and computationally more
useful nonlinear sequence transformations are not mentioned at all.
Apparently, many mathematicians still prefer to work on the theoretically
very nice, but computationally at most moderately powerful regular matrix
transformations.

Nevertheless, there are encouraging signs that the situation is changing
for the better, and there are now several books by mathematicians that
describe how nonlinear sequence transformations can be employed
effectively as computational tools. For example, the most recent (third)
edition of the book \emph{Numerical Recipes}
\cite{Press/Teukolsky/Vetterling/Flannery/2007} now also discusses
nonlinear sequence transformations (for a discussion of the topics
treated there and for further details, see \cite{Weniger/2007d}).

I can also recommend a recent book by Bornemann, Laurie, Wagon, and
Waldvogel \cite{Bornemann/Laurie/Wagon/Waldvogel/2004} on extreme digit
hunting in the context of some challenging problems of numerical
analysis. For this extreme digit hunting, the authors also use sequence
transformations, whose basic theory is described compactly in their
Appendix A. This Appendix is too short to provide a reasonably complete
and balanced presentation of sequence transformation, but I think that a
novice can benefit considerably from reading it. I also like the
extremely pragmatic approach of the authors of this book, which is very
uncommon among mathematicians. Probably, this is due to the fact that the
authors are not primarily interested in the mathematical theory of
sequence transformations: They only wanted to apply sequence
transformations as computational tools in order to obtain more precise
results at tolerable computational costs.

Then, there is a very recent book by Gil, Segura, and Temme
\cite{Gil/Segura/Temme/2007} on the evaluation of special functions. It
discusses in addition to various other computational techniques also
Pad\'{e} approximants, continued fractions, and nonlinear sequence
transformations which all facilitate the evaluation of (power) series
representations for special functions.

My own research shows that nonlinear sequence transformations can be
extremely useful in a large variety of different contexts. I applied them
successfully in such diverse fields as the evaluation of molecular
multicenter integrals of exponentially decaying functions
\cite{Grotendorst/Weniger/Steinborn/1986,Homeier/Weniger/1995,%
  Steinborn/Weniger/1990,Weniger/Grotendorst/Steinborn/1986a,%
  Weniger/Steinborn/1988}, the evaluation of special functions and
related objects \cite{Jentschura/Gies/Valluri/Lamm/Weniger/2002,%
  Jentschura/Mohr/Soff/Weniger/1999,Weniger/1989,Weniger/1990,%
  Weniger/1994a,Weniger/1994b,Weniger/1996d,Weniger/2001,Weniger/2003,%
  Weniger/Cizek/1990,Weniger/Steinborn/1989a}, the summation of strongly
divergent quantum mechanical perturbation expansions
\cite{Bender/Weniger/2001,Cizek/Vinette/Weniger/1991,%
  Cizek/Vinette/Weniger/1993a,Cizek/Weniger/Bracken/Spirko/1996,%
  Jentschura/Becher/Weniger/Soff/2000,Jentschura/Weniger/Soff/2000,%
  Weniger/1990,Weniger/1992,Weniger/1994b,Weniger/1996a,%
  Weniger/1996c,Weniger/1996e,Weniger/1997,Weniger/2001,%
  Weniger/Cizek/Vinette/1991,Weniger/Cizek/Vinette/1993}, the prediction
of unknown perturbation series coefficients
\cite{Bender/Weniger/2001,Jentschura/Becher/Weniger/Soff/2000,%
  Jentschura/Weniger/Soff/2000,Weniger/1997,Weniger/2000b}, and the
extrapolation of quantum chemical crystal orbital and cluster electronic
structure calculations for oligomers to their infinite chain limits of
stereoregular \emph{quasi}-onedimensional organic polymers
\cite{Cioslowski/Weniger/1993,Weniger/Kirtman/2003,Weniger/Liegener/1990}.

In view of all these examples, it is probably justified to claim that
anybody involved in computational work should have at least some basic
knowledge about the power and also about the shortcomings and limitations
of nonlinear sequence transformations.

Pad\'{e} approximants $[m/n]_{f} (z)$ can be viewed to be a special class
of nonlinear sequence transformation since they convert the partial sums
$f_{n} (z) = \sum_{k=0}^{n} \gamma_{k} z^{k}$ of a (formal) power series
for some function $f (z)$ to a doubly indexed sequence of rational
functions. As documented by the long list of successful applications in
the monograph by Baker and Graves-Morris \cite{Baker/Graves-Morris/1996},
Pad\'{e} approximants are now almost routinely used in theoretical
physics and in applied mathematics to overcome problems with slowly
convergent or divergent power series. It is, however, not nearly so well
known among nonspecialists that alternative sequence transformations can
at least for certain computational problems be much more effective than
Pad\'{e} approximants.

The probably best known sequence transformation is Wynn's epsilon
algorithm \cite{Wynn/1956a}, which is defined by the following nonlinear
recursive scheme:
\begin{subequations}
  \label{eps_al}
  \begin{align}
    \label{eps_al_a}
    \epsilon_{-1}^{(n)} & \; = \; 0 \, ,
    \qquad \epsilon_0^{(n)} \, = \, s_n \, ,
    \qquad  n \in \mathbb{N}_0 \, , \\
    \label{eps_al_b}
    \epsilon_{k+1}^{(n)} & \; = \; \epsilon_{k-1}^{(n+1)} \, + \,
    \frac{1}{\epsilon_{k}^{(n+1)} - \epsilon_{k}^{(n)}} \, ,
    \qquad k, n \in \mathbb{N}_0 \, .
  \end{align}
\end{subequations}
The elements $\epsilon_{2k}^{(n)}$ with \emph{even} subscripts provide
approximations to the (generalized) limit $s$ of the sequence $\{ s_n
\}_{n=0}^{\infty}$ to be transformed, whereas the elements
$\epsilon_{2k+1}^{(n)}$ with \emph{odd} subscripts are only auxiliary
quantities which diverge if the whole process converges. A compact
FORTRAN 77 program for the epsilon algorithm as well as the underlying
computational algorithm is described in \cite[Section 4.3]{Weniger/1989}.
In \cite[p.\ 213]{Press/Teukolsky/Vetterling/Flannery/2007}, a
translation of this FORTRAN 77 program to \texttt{C} can be found.

If the elements of the input sequence $\{ s_n \}_{n=0}^{\infty}$ are the
partial sums $f_{n} (z) = \sum_{k=0}^{n} \gamma_{k} z^{k}$ of the
(formal) power series for some function $f (z)$, then the epsilon
algorithm with $\epsilon_{0}^{(n)} = f_{n} (z)$ produces Pad\'{e}
approximants to $f (z)$:
\begin{equation}
  \label{Eps_Pade}
  \epsilon_{2 k}^{(n)} \; = \; [ n + k / k ]_f (z) \, ,
  \qquad k, n \in \mathbb{N}_{0} \, .
\end{equation}
But Wynn's epsilon algorithm is not restricted to input data that are the
partial sums of a (formal) power series. Therefore, it is more general
and more widely applicable than Pad\'{e} approximants. Moreover, the
epsilon algorithm can be generalized to cover for example vector
sequences. A recent review can be found in
\cite{Graves-Morris/Roberts/Salam/2000a}.

Since the epsilon algorithm can be used for the computation of Pad\'e
approximants, it is discussed in books on Pad\'e approximants such as the
one by Baker and Graves-Morris \cite{Baker/Graves-Morris/1996}, but there
is also an extensive literature dealing directly with it. On p.\ 120 of
Wimps book \cite{Wimp/1981} it is mentioned that over 50 articles on the
epsilon algorithm were published by Wynn alone, and at least 30 articles
by Brezinski. As a fairly complete source of references on the epsilon
algorithm, Wimp recommends Brezinski's first book \cite{Brezinski/1977}.
However, this book was published in 1977, and since then many more
articles on the theory or on applications of Wynn's epsilon algorithm
have appeared. Thus, any attempt of providing a reasonably complete
bibliography would be beyond the scope of this article.

In a convergence acceleration or summation process, it is usually a good
idea to try to use the available information as effectively as possible
(possible exceptions to this rule are discussed in \cite{Weniger/2001}).
Let us assume that a finite subset $\bigl\{ s_{0}, s_{1}, \dots, s_{m}
\bigr\}$ of sequence elements is available. Then, those elements
$\epsilon_{2k}^{(n)}$ of the epsilon table produced by (\ref{eps_al})
should be chosen as approximations to the limit of the imput sequence
that have the highest possible subscript or transformation order. Thus,
if $m \in \mathbb{N}_{0}$ is even, $m=2\mu$, I use as approximation to
the limit of the input sequence the transformation \cite[Eq.\
(4.3-4)]{Weniger/1989}
\begin{equation}
  \label{EpsAl_ApprLim_even}
  \bigl\{ s_{0}, s_{1}, \dots, s_{2\mu} \bigr\} \; \to \;
  \epsilon_{2 \mu}^{(0)} \, ,
\end{equation}
and if $m \in \mathbb{N}_{0}$ is odd, $m=2\mu+1$, I use the
transformation \cite[Eq.\ (4.3-5)]{Weniger/1989}
\begin{equation}
  \label{EpsAl_ApprLim_odd}
  \bigl\{ s_{1}, s_{2}, \dots, s_{2\mu+1} \bigr\} \; \to \;
  \epsilon_{2 \mu}^{(1)} \, .
\end{equation}
With the help of the notation $\Ent {x}$ for the integral part of $x$,
which is the largest integer $\nu$ satisfying $\nu \le x$, these two
relationships can be combined into a single equation yielding \cite[Eq.\
(4.3-6)]{Weniger/1989}:
\begin{equation}
  \label{EpsAl_ApprLim_gen}
  \bigl\{ s_{m - 2 \Ent {m/2}}, s_{m - 2 \Ent {m/2} + 1}, \dots ,
  s_{m} \bigr\} \; \to \;
  \epsilon_{2 \Ent {m/2}}^{(m - 2 \Ent {m/2})} \, .
\end{equation}

Wynn's epsilon algorithm is an example of a sequence transformation that
uses as input data only the elements of the sequence to be transformed.
However, in some cases structural information on the dependence of the
sequence elements $s_{n}$ on the index $n$ is available. For example, it
is well known that the truncation error of a convergent series with
strictly alternating and monotonously decreasing terms is bounded in
magnitude by the first term not included in the partial sum and that it
possesses the same sign as this term (see for instance \cite[p.\
259]{Knopp/1964}). The first term neglected is also the best simple
estimate for the truncation error of a strictly alternating
nonterminating and thus diverging hypergeometric series ${}_2 F_0
(\alpha, \beta; - z)$ with $\alpha, \beta, z > 0$ \cite[Theorem
5.12-5]{Carlson/1977}. Such an information on the index dependence of the
truncation errors should be extremely helpful in a convergence
acceleration or summation process, but a sequence transformations like
Wynn's epsilon algorithm cannot benefit from it.

A convenient way of incorporating such an information into the
transformation process consists in the use of remainder estimates $\{
\omega_n \}_{n=0}^{\infty}$. Because of the additional information
contained in the remainder estimates, sequence transformations of that
kind are potentially very powerful as well as very versatile.

The best-known example of such a sequence transformation is Levin's
transformation \cite{Levin/1973}, which is generally considered to be a
very powerful as well as very versatile sequence transformation (see for
example \cite{Brezinski/RedivoZaglia/1991a,Homeier/2000a,%
  Smith/Ford/1979,Smith/Ford/1982,Weniger/1989,Weniger/2004} and
references therein):
\begin{equation}
  \label{GenLevTr}
  \mathcal{L}_{k}^{(n)} (\beta, s_n, \omega_n) \; = \; \frac
{\displaystyle
\sum_{j=0}^{k} \, (-1)^{j} \, {\binom{k}{j}} \,
\frac {(\beta+n+j)^{k-1}} {(\beta+n+k)^{k-1}} \,
\frac {s_{n+j}} {\omega_{n+j}} }
{\displaystyle
\sum_{j=0}^{k} \, (-1)^{j} \, {\binom{k}{j}} \,
\frac {(\beta+n+j)^{k-1}} {(\beta+n+k)^{k-1}} \,
\frac {1} {\omega_{n+j}} }
\, , \qquad k, n \in \mathbb{N}_0 \, .
\end{equation}
Here, $\beta > 0$ is a shift parameter. The most obvious choice is $\beta
= 1$, which is exclusively used in this article.

The numerator and denominator sums of $\mathcal{L}_{k}^{(n)} (\beta, s_n,
\omega_n)$ can also be computed recursively (\cite[Eq.\ (7.2-8) -
(7.2-10)]{Weniger/1989} or in \cite[Eq.\ (3.11)]{Weniger/2004}):
\begin{subequations}
  \label{RecSchemeLevTr}
  \begin{align}
    \label{RecSchemeLevTr_a}
    L_{0}^{(n)} & \; = \; u_n \, ,
    \quad n \in \mathbb{N}_0 \, , \\
    \label{RecSchemeLevTr_b}
    L_{k+1}^{(n)} & \; = \; L_k^{(n+1)} \, - \,
    \frac {(\beta + n) (\beta+n+k)^{k-1}}{(\beta+n+k+1)^k} \,
    L_k^{(n)} \, , \quad k, n \in \mathbb{N}_0 \, ,
  \end{align}
\end{subequations}
By choosing in (\ref{RecSchemeLevTr_a}) either $u_n = s_n/\omega_n$ or
$u_n = 1/\omega_n$, we obtain the numerator and denominator sums of
Levin's transformation (\ref{GenLevTr}).

As discussed in more detail in \cite{Weniger/1989,Weniger/2004}, Levin's
transformation is based on the implicit assumption that the ratio $[s_n -
s]/\omega_n$ can be expressed as a power series in $1/(n+\beta)$.  A
different class of sequence transformations can be derived by assuming
that $[s_n - s]/\omega_n$ can be expressed as a so-called factorial
series, yielding \cite[Eq.\ (8.2-7)]{Weniger/1989}
\begin{equation}
  \label{GenWenTr}
  \mathcal{S}_{k}^{(n)} (\beta, s_n, \omega_n) \; = \; \frac
  {\displaystyle
    \sum_{j=0}^{k} \, (-1)^{j} \, {\binom{k}{j}} \,
    \frac {(\beta+n+j)_{k-1}} {(\beta+n+k)_{k-1}} \,
    \frac {s_{n+j}} {\omega_{n+j}} }
  {\displaystyle
    \sum_{j=0}^{k} \, (-1)^{j} \, {\binom{k}{j}} \,
    \frac {(\beta+n+j)_{k-1}} {(\beta+n+k)_{k-1}} \,
    \frac {1} {\omega_{n+j}} } \, ,
  \qquad k, n \in \mathbb{N}_0 \, .
\end{equation}
As in the case of Levin's transformation. $\beta > 0$ is a shift
parameter, and again, only $\beta=1$ is considered in this article.
Formally, we obtain $\mathcal{S}_{k}^{(n)} (\beta, s_n, \omega_n)$ from
$\mathcal{L}_{k}^{(n)} (\beta, s_n, \omega_n)$ if we replace in
(\ref{GenLevTr}) the powers $(\beta+n+j)^{k-1}$ by Pochhammer symbols
$(\beta+n+j)_{k-1}$.

The numerator and denominator sums of $\mathcal{S}_{k}^{(n)} (\beta, s_n,
\omega_n)$ can also be computed recursively (\cite[Eq.\ (7.2-8) -
(8.3-7)]{Weniger/1989} or in \cite[Eq.\ (3.12)]{Weniger/2004}):
\begin{subequations}
  \label{RecSchemeWenTr}
  \begin{align}
    \label{RecSchemeWenTr_a}
    S_{0}^{(n)} & \; = \; u_n \, ,
    \quad n \in \mathbb{N}_0 \, , \\
    \label{RecSchemeWenTr_b}
    S_{k+1}^{(n)} & \; = \; S_k^{(n+1)} \, - \,
    \frac {(\beta+n+k-1) (\beta + n + k)} {(\beta+n+2k-1) (\beta+n+2k)}
    \, S_k^{(n)} \, , \quad k, n \in \mathbb{N}_0 \, ,
  \end{align}
\end{subequations}
The initial values $u_n = s_n/\omega_n$ produce the numerators of the
transformation (\ref{GenWenTr}), and the initial values $u_n =
1/\omega_n$ yield the denominators.

Both Levin's sequence transformation (\ref{GenLevTr}) as well as the
related transformation (\ref{GenWenTr}) utilize the information contained
in explicit remainder estimates $\{ \omega_n \}_{n=0}^{\infty}$ which
should be chosen in such a way that the ratio $[s_{n}-s]/\omega_{n}$
becomes a smooth function of $n$ that can be annihilated effectively by
weighted finite difference operators (see for example \cite[Sections II
and IV]{Weniger/2004} and references therein). Accordingly, the choice of
the remainder estimates is of utmost importance for the success or
failure of a convergence acceleration or summation process involving
Levin-type transformations. In this respect, it may be interesting to
note that a symbolic approach for the construction of asymptotic
estimates to the truncation errors of series representations for special
functions was recently developed in \cite{Weniger/2007a}.

On the basis of purely heuristic arguments Levin \cite{Levin/1973} had
suggested some simple remainder estimates which according to experience
work remarkably well in a large variety of cases and which can also be
used in the case of the sequence transformation (\ref{GenWenTr}).  But
the best \emph{simple} estimate for the truncation error of a strictly
alternating convergent series is the first term not included in the
partial sum \cite[p.\ 259]{Knopp/1964}.  Moreover, the first term
neglected is also an estimate of the truncation error of a divergent
hypergeometric series ${}_2 F_0 (a, b, - z)$ with $a, b, z > 0$
\cite[Theorem 5.12-5]{Carlson/1977}. Accordingly, Smith and Ford
\cite{Smith/Ford/1979} proposed for sequences of partial sums of
alternating series the following remainder estimate:
\begin{equation}
  \label{dRemEst}
  \omega_n \; = \; \Delta s_n \, .
\end{equation}
The use of this remainder estimate in \cite{Levin/1973} and
(\ref{GenWenTr}) yields the following variants of the sequence
transformations $\mathcal{L}_{k}^{(n)} (\beta, s_n, \omega_n)$ and
$\mathcal{S}_{k}^{(n)} (\beta, s_n, \omega_n)$:
\begin{align}
  \label{dLevTr}
  d_{k}^{(n)} (\beta, s_n) & \; = \;
  \mathcal{L}_{k}^{(n)} (\beta, s_n, \Delta s_n) \, ,
  \\
  \label{dWenTr}
  \delta_{k}^{(n)} (\beta, s_n) & \; = \;
  \mathcal{S}_{k}^{(n)} (\beta, s_n, \Delta s_n) \, .
\end{align}
If the elements of the input sequence $\{ s_n \}_{n=0}^{\infty}$ are the
partial sums $f_{n} (z) = \sum_{k=0}^{n} \gamma_{k} z^{k}$ of the
(formal) power series for some function $f (z)$, then $d_{k}^{(n)} \bigl(
\beta, f_{n} (z) \bigr)$ and $\delta_{k}^{(n)} \bigl( \beta, f_{n} (z)
\bigr)$ are ratios of two polynomials in $z$ of degrees $k+n$ and $k$,
respectively (for a detailed discussion of these rational approximants,
see \cite[Section VI]{Weniger/2004}).

The transformations $d_{k}^{(n)} (\beta, s_n)$ and in particular also
${\delta}_{k}^{(n)} (\beta, s_n)$ were found to be remarkably powerful
summation techniques for divergent series with strictly alternating
terms. Numerous successful applications of Levin's sequence
transformation (\ref{GenLevTr}) and of the related transformation
(\ref{GenWenTr}) and their variants in convergence acceleration and
summation processes are described in \cite[pp.\ 1210 and
1225]{Weniger/2004}. It may be interesting to note that recently
${\delta}_{k}^{(n)} (\beta, s_n)$ has been used quite a lot in optics
\cite{Borghi/2007,Borghi/2008,Borghi/Alonso/2007,Borghi/Santarsiero/2003}.

In the case of the transformations (\ref{dLevTr}) and (\ref{dWenTr}),
the approximation to the limit with the highest transformation order is
given by
\begin{equation}
  \label{LevinTypeApprLim}
\{ s_0, s_1, \ldots ,s_{m+1} \} \; \to \;
\Xi_{m}^{(0)} (\beta, s_0) \, ,
\end{equation}
where $\Xi_{k}^{(n)} (\beta, s_n)$ stands for either $d_{k}^{(n)}
(\beta, s_n)$ or ${\delta}_{k}^{(n)} (\beta, s_n)$.
\end{appendix}

\addcontentsline{toc}{section}{Bibliography}

\providecommand{\SortNoop}[1]{} \providecommand{\OneLetter}[1]{#1}
  \providecommand{\SwapArgs}[2]{#2#1}

\end{document}